\documentclass[a4paper]{amsart}

\usepackage{amsmath,amssymb,amsthm,mathtools,bm}
\usepackage{hyperref}
\usepackage[abs]{overpic}

\newtheorem{theorem}{Theorem}[section]
\newtheorem{lemma}[theorem]{Lemma}
\newtheorem{proposition}[theorem]{Proposition}
\newtheorem{corollary}[theorem]{Corollary}

\theoremstyle{definition}

\newtheorem{example}[theorem]{Example}
\newtheorem{remark}[theorem]{Remark}

\theoremstyle{remark}
\newtheorem*{acknowledgments}{Acknowledgments}
\newtheorem*{declaration}{Declaration of AI use}

\newcommand{\As}{\operatorname{As}}
\newcommand{\Conj}{\operatorname{Conj}}
\newcommand{\Supp}{\operatorname{Supp}}
\newcommand{\id}{\operatorname{id}}
\newcommand{\Z}{\mathbb Z}
\newcommand{\Sym}{\mathfrak S}

\makeatletter

\@addtoreset{figure}{section}
\makeatother

\makeatletter

\@addtoreset{table}{section}
\makeatother

\makeatletter
  
  \@addtoreset{equation}{section}
\makeatother

\title[Hurwitz equivalence in quandles of order at most four]{Hurwitz equivalence in quandles of order at most four}

\author{Takuto Hosoo}
\address{Data Analytics Labo Co., 6F, Shiodome Eastside Building, 5-4-18 Tsukiji, Chuo-ku, Tokyo 104-0045, Japan}
\email{243s017s@gsuite.kobe-u.ac.jp}

\author[Kodai Wada]{Kodai Wada}
\address{Department of Mathematics, Kobe University, 
1-1 Rokkodai-cho, Nada-ku, Kobe 657-8501, Japan}
\email{wada@math.kobe-u.ac.jp}

\makeatletter
\@namedef{subjclassname@2020}{%
  \textup{2020} Mathematics Subject Classification}
\makeatother


\begin{document}

\begin{abstract}
Let $X$ be a quandle and let $m\geq 2$ be an integer. 
The $m$-braid group acts on the $m$-fold Cartesian power $X^m$ by the Hurwitz action. 
For an $m$-tuple $\bm{x}\in X^m$, we introduce two invariants of the Hurwitz action: 
the product $\Pi(\bm{x})$ of its entries in the associated group $\As(X)$ and the orbit multiplicities $\mu_w(\bm{x})$ $(w\in X)$. 
We prove that these two invariants completely classify Hurwitz equivalence classes for all quandles of order at most four. 
We also establish a complete classification of Hurwitz equivalence classes 
for a class of amalgamated sums of trivial quandles and a reduction theorem for disjoint sums with trivial quandles. 
Finally, we study the corresponding action of the virtual braid group. 
For this action, the orbit multiplicities alone give a complete classification for every quandle of order at most four.
\end{abstract}

\maketitle

\section{Introduction}

Let $m\geq 2$ be an integer, and
let $B_m$ denote the braid group on $m$ strands
with standard generators $\sigma_1,\ldots,\sigma_{m-1}$.
For a group $G$, the classical Hurwitz action of $B_m$ on the $m$-fold Cartesian power $G^m$ is defined by
\[
(g_1,\ldots,g_i,g_{i+1},\ldots,g_m)\cdot\sigma_i
=
(g_1,\ldots,g_{i-1},
g_{i+1},g_{i+1}^{-1}g_i g_{i+1},
g_{i+2},\ldots,g_m)\in G^m.
\]
This action preserves the product $g_1\cdots g_m$
and hence acts on the set of factorizations
of a fixed element of $G$.
Two factorizations are said to be Hurwitz equivalent
if they belong to the same orbit under this action.
Hurwitz equivalence of such factorizations
arises naturally in braid monodromy, 
where changes of geometric bases are described by the Hurwitz action of the braid group. 
Hurwitz equivalence of braid monodromy factorizations is 
related to the topology of plane curves and algebraic
surfaces~\cite{KT}, 
and also appears naturally in the study of braided surfaces~\cite{Kam96,Kam02}. 

The problem of classifying Hurwitz equivalence classes has been extensively studied in combinatorial group theory. 
Complete classifications have been obtained in several settings, 
including certain factorizations into transpositions
in symmetric groups~\cite{BT},
tuples of elements in dicyclic and semidihedral groups~\cite{Sia},
and tuples of elements in dihedral groups~\cite{Ber,Sia}.

The Hurwitz action admits a natural formulation for arbitrary quandles.
Recall that a quandle is a nonempty set $X$ equipped with a binary operation
$*$ modeled on conjugation in a group~\cite{Joy,Mat}.
Quandles arise naturally in knot theory, 
where their axioms reflect the algebraic structure underlying
Reidemeister moves.
For a quandle $X$, the braid group $B_m$ acts on the $m$-fold Cartesian power $X^m$ by 
\[
(x_1,\ldots,x_i,x_{i+1},\ldots,x_m)\cdot\sigma_i
=
(x_1,\ldots,x_{i-1},
x_{i+1},x_i*x_{i+1},
x_{i+2},\ldots,x_m)\in X^m.
\]

Compared with the group case, however, much less is known about the orbit structure of the Hurwitz action on quandles. 
Among the cases studied in the group setting, 
the results for dihedral groups immediately yield classifications for finite dihedral quandles: 
the set of reflections in the dihedral group $D_{2k}$ of order $2k$ forms a quandle isomorphic to the dihedral quandle $R_k$ of order $k$. 
Thus the results of Berger~\cite{Ber} and Sia~\cite{Sia} give complete classifications for finite dihedral quandles. 
More recently, Nakamura, Nakanishi, Satoh, and the second author
\cite{NNSW} obtained a complete classification for the infinite
dihedral quandle $R_\infty$.

In this paper, we study Hurwitz equivalence for arbitrary quandles
and obtain a complete classification of Hurwitz equivalence classes for all quandles of order at most four.
Our approach is based on two invariants that are defined for every
quandle $X$.
For $\bm{x}=(x_1,\ldots,x_m)\in X^m$, let $\As(X)$ be the
associated group of $X$, and set
\[
\Pi(\bm{x})=x_1\cdots x_m\in\As(X), 
\]
where each $x_i$ also denotes its image under the natural map $X\to\As(X)$. 
For each $w\in X$, let $O_{\bm{x}}(w)$ be the orbit of $w$
under the operations $*x_i$ and $*^{-1}x_i$ $(1\leq i\leq m)$.
We define the \emph{orbit multiplicity} of $w$ in $\bm{x}$ by
\[
\mu_w(\bm{x})
=\#\{i\mid x_i\in O_{\bm{x}}(w)\}.
\]
These invariants are quandle analogues of those used by Berger~\cite{Ber} in the group setting. 
We show that $\Pi(\bm{x})$ and $\mu_w(\bm{x})$ are preserved
under Hurwitz equivalence for every quandle~$X$.

Up to isomorphism, there are twelve quandles of order at most four, denoted by 
\[
X_{1,1},\ X_{2,1},\ X_{3,1},\ X_{3,2},\ X_{3,3},
\ X_{4,1},\ldots,X_{4,7}.
\]
This collection includes several different types of examples: 
$X_{3,3}$ and $X_{4,6}$ are the dihedral quandles of orders three and four, respectively, 
while $X_{4,7}$ is the tetrahedral quandle. 
Our first main result shows that the two invariants introduced above form a complete set of invariants for distinguishing Hurwitz equivalence classes in all these cases.

\begin{theorem}\label{thm:main-classical}
Let $X$ be an arbitrary quandle of order at most four, and
let $\bm{x},\bm{y}\in X^m$.
Then the following are equivalent.
\begin{enumerate}
\item
There exists a braid $\beta\in B_m$ such that
$\bm{x}\cdot\beta=\bm{y}$.
\item
$\Pi(\bm{x})=\Pi(\bm{y})$ and
$\mu_w(\bm{x})=\mu_w(\bm{y})$ for every $w\in X$.
\end{enumerate}
\end{theorem}

In particular, Theorem~\ref{thm:main-classical} includes the dihedral quandles of orders three and four, 
and thus recovers, by a different approach, the corresponding classifications that follow from the results of Berger~\cite{Ber} and Sia~\cite{Sia}.

Beyond the case of quandles of order at most four, 
we establish a complete classification for a class of amalgamated sums 
and a reduction theorem for disjoint sums with trivial quandles. 
First, we consider a class of \emph{amalgamated sums}
arising from two finite trivial quandles $Y$ and $Z$.
Such a quandle $X=Y\sqcup Z$ is determined by
a family of permutations $\theta_z$ of $Y$ associated with $z\in Z$ 
such that $\theta_z\theta_{z'}=\theta_{z'}\theta_z$ for all $z,z'\in Z$.
The operation is defined by $y*z=\theta_z(y)$
for $y\in Y$ and $z\in Z$,
and is trivial in all other cases.
For every quandle in this class, the orbit multiplicities alone 
completely determine Hurwitz equivalence classes 
(Theorem~\ref{thm:amalgamated-sum-classification}).
In particular, this gives a uniform classification for 
$X_{3,2}$, $X_{4,2}$, $X_{4,3}$, and $X_{4,4}$. 

Second, let $X=Q\sqcup Z$ be the \emph{disjoint sum}
of an arbitrary quandle $Q$ and a finite trivial quandle $Z$,
with trivial operations between the two components.
For a tuple $\bm{x}\in X^m$, write $\widehat{\bm{x}}$
for the tuple obtained by deleting all entries belonging to $Z$.
We show that two tuples $\bm{x},\bm{y}\in X^m$
are Hurwitz equivalent in $X$ if and only if
$\mu_z(\bm{x})=\mu_z(\bm{y})$ for every $z\in Z$
and $\widehat{\bm{x}}$ and $\widehat{\bm{y}}$
are Hurwitz equivalent in $Q$.
Thus the classification for $X$ reduces to that for $Q$,
together with the multiplicities of the entries in $Z$
(Theorem~\ref{thm:disjoint-sum-classification}).
We apply this reduction to $X_{4,5}$, which is the disjoint
sum of $X_{3,3}$ and a one-element trivial quandle.

Although the product $\Pi(\bm{x})$ is conceptually natural, it is not always easy to determine whether 
$\Pi(\bm{x})=\Pi(\bm{y})$, since these products take values in
$\As(X)$.
To make this invariant more tractable, for a quandle homomorphism
$\phi:X\to\Conj(G)$ into the conjugation quandle of a group $G$, we
introduce
\[
\Pi_\phi(\bm{x})
=
\phi(x_1)\cdots\phi(x_m)\in G.
\]
By the universal property of $\As(X)$,
the homomorphism $\phi$ induces a unique group homomorphism
from $\As(X)$ to $G$, under which
$\Pi(\bm{x})$ maps to $\Pi_\phi(\bm{x})$.
Suitable choices of $\phi$ allow us to distinguish
the normal forms arising in the classification.

This point of view also connects our invariants with the previously
known classification of the infinite dihedral quandle.
In Section~\ref{sec:invariants}, we clarify the relationship between the invariants
$\Delta(\bm{x})$, $d(\bm{x})$, and $M(\bm{x})$ used in~\cite{NNSW}
and our invariants $\Pi_\phi(\bm{x})$ and $\mu_w(\bm{x})$ 
(Theorem~\ref{thm:Rinfty-invariants} and Corollary~\ref{cor:Rinfty-invariants}).

We also consider the Hurwitz action of the virtual braid group $VB_m$ on $X^m$. 
In addition to the braid generators $\sigma_i$, 
the group $VB_m$ has generators $\tau_i$, 
where each $\tau_i$ acts on $X^m$ by interchanging the $i$-th and $(i+1)$-st entries. 
Although the product $\Pi(\bm{x})$ is no longer preserved in general under these transpositions, 
the orbit multiplicities remain unchanged. 
Our second main result shows that, for quandles of order at most four,
the orbit multiplicities alone completely classify virtual Hurwitz equivalence classes.

\begin{theorem}\label{thm:main-virtual}
Let $X$ be an arbitrary quandle of order at most four, and
let $\bm{x},\bm{y}\in X^m$.
Then the following are equivalent.
\begin{enumerate}
\item
There exists a virtual braid $\beta\in VB_m$ such that
$\bm{x}\cdot\beta=\bm{y}$.
\item
$\mu_w(\bm{x})=\mu_w(\bm{y})$ for every $w\in X$.
\end{enumerate}
\end{theorem}

We establish this result by constructing explicit normal forms
for virtual Hurwitz equivalence classes.

The paper is organized as follows.
In Section~\ref{sec:quandles}, we recall basic facts about
quandles and associated groups, and list the twelve quandles
of order at most four.
In Section~\ref{sec:invariants}, we define the Hurwitz action
and introduce the invariants $\Pi(\bm{x})$ and
$\mu_w(\bm{x})$, as well as $\Pi_\phi(\bm{x})$.
We also compare these invariants with the invariants used for $R_\infty$ in~\cite{NNSW}.
In Section~\ref{sec:sum-constructions}, we classify Hurwitz equivalence classes for trivial quandles and a class of amalgamated sums, 
and establish a reduction theorem for disjoint sums with trivial quandles.
These results provide a complete classification of Hurwitz equivalence classes for the eight quandles 
$X_{1,1}$, $X_{2,1}$, $X_{3,1}$, $X_{3,2}$, $X_{4,1}$, $X_{4,2}$, $X_{4,3}$, and $X_{4,4}$. 
In Section~\ref{sec:dihedral-order3}, we classify Hurwitz equivalence classes for $X_{3,3}$ and apply the reduction theorem to $X_{4,5}$. 
Sections~\ref{sec:dihedral-order4} and \ref{sec:tetrahedral} are devoted to the dihedral quandle $X_{4,6}$ and the tetrahedral quandle $X_{4,7}$, respectively.
Finally, in Section~\ref{sec:virtual}, we study the Hurwitz action of the 
virtual braid group, construct explicit normal
forms, and use them to prove
Theorem~\ref{thm:main-virtual}.

\begin{acknowledgments}
The authors would like to express their sincere gratitude to Katsumi Ishikawa for bringing~\cite{Cla} to their attention,
thereby enabling them to complete the proof of
Theorem~\ref{thm:X47}.
They would also like to thank Shin Satoh for many fruitful discussions and for his encouragement throughout this work.
This work was supported by JSPS KAKENHI Grant Number JP23K12973.
\end{acknowledgments}

\begin{declaration}
ChatGPT (OpenAI) was used to assist with language editing, 
the exposition of the manuscript, and the review of mathematical arguments and calculations. 
All mathematical statements and proofs were independently checked by the authors, who take full responsibility for the final content.
\end{declaration}

\section{Preliminaries on quandles and associated groups}
\label{sec:quandles}

A \emph{quandle} is a pair $(X,*)$ consisting of a nonempty set $X$ and a binary operation $*:X\times X\to X$ satisfying the following axioms (see~\cite{Joy,Mat}):
\begin{enumerate}
\item 
$x*x=x$ for all $x\in X$; 
\item 
for all $x,y\in X$, there exists a unique $z\in X$ such that
$z*y=x$; 
\item 
$(x*y)*z=(x*z)*(y*z)$ for all $x,y,z\in X$. 
\end{enumerate}

The \emph{order} of $X$ is its cardinality. 
The second axiom enables us to define the \emph{dual operation $*^{-1}$} 
as the unique binary operation $*^{-1}:X\times X\to X$ satisfying 
\[
(x*y)*^{-1}y=(x*^{-1}y)*y=x.
\]

A quandle $(A,*_A)$ is a \emph{subquandle} of a quandle $(X,*)$ 
if $A$ is a subset of $X$ and the operation $*_A$ is the restriction of $*$. 
We also denote $*_A$ by $*$.

\begin{example}\label{ex:quandles}
(i) 
Every nonempty set $X$ becomes a quandle under 
\[
x*y=x*^{-1}y=x \quad (x,y\in X).
\]
This quandle is called a \emph{trivial quandle}. 

(ii) 
Let $X=\Z/k\Z$ be the cyclic group of order $k$~$(k\geq3)$.
Then $X$ becomes a quandle under 
\[
x*y=x*^{-1}y=2y-x \quad (x,y\in X).
\]
This quandle is called the \emph{dihedral quandle of order $k$} and is denoted by $R_k$. 

(iii) 
Let $X=\mathbb{F}_2[t^{\pm1}]/(1+t+t^{2})=\{0,1,t,1+t\}$ 
be the quotient ring of the Laurent polynomial ring $\mathbb{F}_2[t^{\pm1}]$ by the ideal generated by $1+t+t^2$. 
Then $X$ becomes a quandle under 
\[
x*y=tx+(1+t)y \quad (x,y\in X).
\]
This quandle is called the \emph{tetrahedral quandle}. 
The dual operation is given by 
\[
x\mathbin{*^{-1}}y=(1+t)x+ty \quad (x,y\in X).
\]

(iv) 
Every group $G$ becomes a quandle under
\[
x*y=y^{-1}xy \quad (x,y\in G).
\]
This quandle is called the \emph{conjugation quandle} of $G$ and is denoted by $\Conj(G)$. 
The dual operation is given by 
\[
x*^{-1}y=yxy^{-1} \quad (x,y\in G).
\]

Let $H$ be a subset of $G$ that is invariant under all conjugations by elements of $G$. 
Then $(H,*)$ is a subquandle of $\Conj(G)$ and 
is denoted by $\Conj(H)$. 
\end{example}

Let $X$ and $Y$ be quandles. 
A map $\phi:X\to Y$ is called a \emph{quandle homomorphism} if for all $x,y\in X$, 
\[
\phi(x*y)=\phi(x)*\phi(y). 
\]
A bijective quandle homomorphism $\phi:X\to Y$ is called an \emph{isomorphism}. 
We say that $X$ and $Y$ are \emph{isomorphic} if there exists an isomorphism between $X$ and $Y$. 

\begin{example}\label{ex:dihedral-isomorphism}
(i) Let 
\[
D_{2k}=\langle r,s\mid r^k=s^2=(sr)^2=e\rangle
=\{e,r,\ldots,r^{k-1},s,sr,\ldots,sr^{k-1}\}
\]
be the dihedral group of order $2k$ with identity element $e$,  
where $r$ is a rotation and $s$ is a reflection of a regular $k$-gon. 
Define a map $\phi_k:R_k\to \Conj(D_{2k})$ by 
\[
\phi_k(x)=sr^x \quad (x\in R_k). 
\]
Since the elements $s,sr,\ldots,sr^{k-1}$ are pairwise distinct in $D_{2k}$, the map $\phi_k$ is injective. 
In $\Conj(D_{2k})$, for $x,y\in R_k$, we have 
\[
sr^x*sr^y=(sr^y)^{-1}sr^xsr^y=sr^{2y-x}.
\]
Hence 
\[
\phi_k(x)*\phi_k(y)=sr^x*sr^y=sr^{2y-x}
=\phi_k(2y-x)=\phi_k(x*y).
\]
Thus $\phi_k$ is a quandle homomorphism. 

(ii) 
Let $H_k$ be the subset of $D_{2k}$ consisting of all reflections, i.e., 
\[
H_k=\{s,sr,\ldots,sr^{k-1}\}. 
\]
Since $H_k$ is closed under conjugation, 
we obtain a subquandle $\Conj(H_k)$ of $\Conj(D_{2k})$, 
see Example~\ref{ex:quandles}(iv). 
The image of $\phi_k$ coincides with $H_k$, 
and hence $\phi_k$ induces an isomorphism $R_k\to\Conj(H_k)$.  
\end{example}

There exist twelve quandles of order at most four up to isomorphism~(see \cite[Appendix]{CKS} and \cite[Section~8.3]{Kam17}), 
denoted by 
\[X_{1,1},\ X_{2,1},\ X_{3,1},\ X_{3,2},\ X_{3,3},\ 
X_{4,1},\ldots,X_{4,7}.
\]  
See Tables~\ref{tab:qdl-order3} and~\ref{tab:qdl-order4}. 
Each quandle is specified by its operation table, 
whose $(i,j)$-entry is $i*j$. 
We remark that $X_{1,1}$, $X_{2,1}$, $X_{3,1}$, and $X_{4,1}$ are trivial quandles. 
The quandles $X_{3,3}$ and $X_{4,6}$ are isomorphic to the dihedral quandles of orders three and four, respectively. 
The quandle $X_{4,7}$ is isomorphic to the tetrahedral quandle. 

\begin{table}[htbp]
\caption{Quandles of order at most three}
\centering

\begin{tabular}{c@{\qquad\qquad}c}
\begin{tabular}{c}
$X_{1,1}$ \\[2mm]
\begin{tabular}{c|c}
$\ast$ & 1 \\ \hline
1 & 1
\end{tabular}
\end{tabular}
&
\begin{tabular}{c}
$X_{2,1}$ \\[2mm]
\begin{tabular}{c|cc}
$\ast$ & 1 & 2 \\ \hline
1 & 1 & 1 \\
2 & 2 & 2
\end{tabular}
\end{tabular}
\end{tabular}

\vspace{5mm}

\begin{tabular}{c@{\qquad}c@{\qquad}c}
\begin{tabular}{c}
$X_{3,1}$ \\[2mm]
\begin{tabular}{c|ccc}
$\ast$ & 1 & 2 & 3 \\ \hline
1 & 1 & 1 & 1 \\
2 & 2 & 2 & 2 \\
3 & 3 & 3 & 3
\end{tabular}
\end{tabular}
&
\begin{tabular}{c}
$X_{3,2}$ \\[2mm]
\begin{tabular}{c|ccc}
$\ast$ & 1 & 2 & 3 \\ \hline
1 & 1 & 1 & 2 \\
2 & 2 & 2 & 1 \\
3 & 3 & 3 & 3
\end{tabular}
\end{tabular}
&
\begin{tabular}{c}
$X_{3,3}$ \\[2mm]
\begin{tabular}{c|ccc}
$\ast$ & 1 & 2 & 3 \\ \hline
1 & 1 & 3 & 2 \\
2 & 3 & 2 & 1 \\
3 & 2 & 1 & 3
\end{tabular}
\end{tabular}
\end{tabular}
\label{tab:qdl-order3}
\end{table}

\begin{table}[htbp]
\caption{Quandles of order four}
\centering

\begin{tabular}{c@{\qquad}c@{\qquad}c}
\begin{tabular}{c}
$X_{4,1}$ \\[2mm]
\begin{tabular}{c|cccc}
$\ast$ & 1 & 2 & 3 & 4 \\ \hline
1 & 1 & 1 & 1 & 1 \\
2 & 2 & 2 & 2 & 2 \\
3 & 3 & 3 & 3 & 3 \\
4 & 4 & 4 & 4 & 4
\end{tabular}
\end{tabular}
&
\begin{tabular}{c}
$X_{4,2}$ \\[2mm]
\begin{tabular}{c|cccc}
$\ast$ & 1 & 2 & 3 & 4 \\ \hline
1 & 1 & 1 & 2 & 1 \\
2 & 2 & 2 & 1 & 2 \\
3 & 3 & 3 & 3 & 3 \\
4 & 4 & 4 & 4 & 4
\end{tabular}
\end{tabular}
&
\begin{tabular}{c}
$X_{4,3}$ \\[2mm]
\begin{tabular}{c|cccc}
$\ast$ & 1 & 2 & 3 & 4 \\ \hline
1 & 1 & 1 & 2 & 2 \\
2 & 2 & 2 & 1 & 1 \\
3 & 3 & 3 & 3 & 3 \\
4 & 4 & 4 & 4 & 4
\end{tabular}
\end{tabular}
\end{tabular}

\vspace{6mm}

\begin{tabular}{c@{\qquad}c@{\qquad}c}
\begin{tabular}{c}
$X_{4,4}$ \\[2mm]
\begin{tabular}{c|cccc}
$\ast$ & 1 & 2 & 3 & 4 \\ \hline
1 & 1 & 1 & 1 & 2 \\
2 & 2 & 2 & 2 & 3 \\
3 & 3 & 3 & 3 & 1 \\
4 & 4 & 4 & 4 & 4
\end{tabular}
\end{tabular}
&
\begin{tabular}{c}
$X_{4,5}$ \\[2mm]
\begin{tabular}{c|cccc}
$\ast$ & 1 & 2 & 3 & 4 \\ \hline
1 & 1 & 3 & 2 & 1 \\
2 & 3 & 2 & 1 & 2 \\
3 & 2 & 1 & 3 & 3 \\
4 & 4 & 4 & 4 & 4
\end{tabular}
\end{tabular}
&
\begin{tabular}{c}
$X_{4,6}$ \\[2mm]
\begin{tabular}{c|cccc}
$\ast$ & 1 & 2 & 3 & 4 \\ \hline
1 & 1 & 1 & 2 & 2 \\
2 & 2 & 2 & 1 & 1 \\
3 & 4 & 4 & 3 & 3 \\
4 & 3 & 3 & 4 & 4
\end{tabular}
\end{tabular}
\end{tabular}

\vspace{6mm}

\begin{tabular}{c}
$X_{4,7}$ \\[2mm]
\begin{tabular}{c|cccc}
$\ast$ & 1 & 2 & 3 & 4 \\ \hline
1 & 1 & 4 & 2 & 3 \\
2 & 3 & 2 & 4 & 1 \\
3 & 4 & 1 & 3 & 2 \\
4 & 2 & 3 & 1 & 4
\end{tabular}
\end{tabular}

\label{tab:qdl-order4}
\end{table}

The \emph{associated group} of a quandle $X$ is defined by 
\[
\As(X)=\langle x\in X \mid x*y=y^{-1}xy\text{ for all }x,y\in X\rangle.
\]
Equivalently,
\[
\As(X)=F(X)/\langle\!\langle \{(x*y)^{-1}y^{-1}xy\mid x,y\in X\}\rangle\!\rangle, 
\]
where $F(X)$ is the free group generated by $X$. 
Let $\eta_X:X\to \As(X)$ denote the natural map;  
that is, $\eta_X$ is the composition of the inclusion map $\iota:X\to F(X)$ and the quotient map $\pi:F(X)\to\As(X)$. 
For simplicity, we write $x$ for $\eta_X(x)$ when no confusion can arise. 

The associated group $\As(X)$ satisfies the following universal property. 

\begin{lemma}[{\cite[Proposition~2.1]{FR}, \cite[Proposition~8.8.4]{Kam17}}]
\label{lem:universality}
For every group $G$ and every quandle homomorphism $\phi:X\to\Conj(G)$, 
there exists a unique group homomorphism $\phi_\#:\As(X)\to G$ such that $\phi_\#\circ\eta_X=\phi$. 
\qed
\end{lemma}

\section{Hurwitz equivalence and invariants}\label{sec:invariants}

\begin{figure}[htbp]
  \vspace{2mm}
  \begin{center}
    \begin{overpic}[]{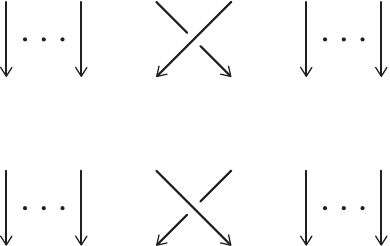}
      \put(-26,96){$\sigma_{i}=$}
      \put(-34,15.7){$\sigma_{i}^{-1}=$}
      \put(-2,124.5){$x_1$}
      \put(31.5,124.5){$x_{i-1}$}
      \put(69,124.5){$x_i$}
      \put(104,124.5){$x_{i+1}$}
      \put(140,124.5){$x_{i+2}$}
      \put(178,124.5){$x_{m}$}
      \put(-2,69){$x_1$}
      \put(31.5,69){$x_{i-1}$}
      \put(62,69){$x_{i+1}$}
      \put(93,69){$x_i*x_{i+1}$}
      \put(140,69){$x_{i+2}$}
      \put(178,69){$x_{m}$}
      \put(-2,43){$x_1$}
      \put(31.5,43){$x_{i-1}$}
      \put(69,43){$x_i$}
      \put(104,43){$x_{i+1}$}
      \put(140,43){$x_{i+2}$}
      \put(178,43){$x_{m}$}
      \put(-2,-11){$x_1$}
      \put(31.5,-11){$x_{i-1}$}
      \put(59,-11){$x_{i+1}*^{-1}x_i$}
      \put(116,-11){$x_i$}
      \put(140,-11){$x_{i+2}$}
      \put(178,-11){$x_{m}$}
    \end{overpic}
  \end{center}
  \vspace{2mm}
  \caption{The Hurwitz action of $B_{m}$ on $X^{m}$}
  \label{fig:Hurwitz-action}
\end{figure}

Let $m\geq 2$ be an integer. 
The braid group $B_m$ on $m$ strands has presentation
\[
B_m=
\left\langle
\sigma_1,\ldots,\sigma_{m-1}
\ \middle|\
\begin{array}{ll}
\sigma_i\sigma_j=\sigma_j\sigma_i,& |i-j|>1,\\
\sigma_i\sigma_{i+1}\sigma_i=\sigma_{i+1}\sigma_i\sigma_{i+1},&
1\leq i\leq m-2
\end{array}
\right\rangle .
\]
Let $X$ be a quandle. 
The \emph{Hurwitz action} of $B_m$ on $X^m$
is defined on the generators by
\[
(x_1,\ldots,x_i,x_{i+1},\ldots,x_m)\cdot \sigma_i
=
(x_1,\ldots,x_{i-1},x_{i+1},x_i*x_{i+1},x_{i+2},\ldots,x_m),
\]
and
\[
(x_1,\ldots,x_i,x_{i+1},\ldots,x_m)\cdot \sigma_i^{-1}
=
(x_1,\ldots,x_{i-1},x_{i+1}*^{-1} x_i,x_i,x_{i+2},\ldots,x_m).
\]
As illustrated in Figure~\ref{fig:Hurwitz-action}, 
this action can be interpreted in terms of
$X$-colorings of braid diagrams.
More precisely, for $\bm{x},\bm{y}\in X^m$ and $\beta\in B_m$, we have $\bm{x}\cdot\beta=\bm{y}$ 
if and only if there exists an $X$-coloring of an oriented diagram of $\beta$ such that the top and bottom endpoints are colored by
$\bm{x}$ and $\bm{y}$, respectively. 
This knot-theoretic interpretation plays an important role in the
proofs of the results below.

Two tuples $\bm{x},\bm{y}\in X^m$ are \emph{Hurwitz equivalent} 
if there exists $\beta\in B_m$ such that 
$\bm{x}\cdot \beta=\bm{y}$. 
We write $\bm{x}\sim \bm{y}$. 
A tuple $\bm{x}\in X^m$ is called \emph{trivial} 
if $\bm{x}=(z,\ldots,z)$ for some $z\in X$. 
We remark that 
if $\bm{x}\sim\bm{y}$ and $\bm{x}$ is trivial, 
then $\bm{x}=\bm{y}$. 
In other words, the orbit of a trivial tuple $\bm{x}$ 
is the singleton $\{\bm{x}\}$.

Let $\bm{x}=(x_1,\ldots,x_m)\in X^m$.
As quandle analogues of Berger's invariants in~\cite{Ber},
we introduce the following two quantities.

\begin{itemize}
\item
The product
\[
\Pi(\bm{x})=x_1\cdots x_m\in\As(X).
\]

\item
For $w\in X$, let $O_{\bm{x}}(w)$ be the smallest
subset of $X$ containing $w$ such that
\[
a*x_i,\ a*^{-1}x_i\in O_{\bm{x}}(w)
\]
for every $a\in O_{\bm{x}}(w)$ and every
$1\le i\le m$.
The \emph{orbit multiplicity} of $w$ in $\bm{x}$ is defined by
\[
\mu_w(\bm{x})
=\#\{i\mid x_i\in O_{\bm{x}}(w)\}.
\]
\end{itemize}

\begin{remark}
For $\bm{x}=(x_1,\ldots,x_m)\in X^m$,
let $\langle\bm{x}\rangle$ denote the subgroup of
$\As(X)$ generated by the entries $x_1,\ldots,x_m$.
Since
\[
x_i*x_{i+1}=x_{i+1}^{-1}x_i x_{i+1}
\]
in $\As(X)$, the subgroup $\langle\bm{x}\rangle$
is preserved under the Hurwitz action.
The quandle operation induces a natural right action
of $\As(X)$ on $X$, determined by
\[
y\cdot z=y*z,
\qquad
y\cdot z^{-1}=y*^{-1}z
\qquad(y,z\in X).
\]
For $w\in X$, the set $O_{\bm{x}}(w)$
can be written as
\[
O_{\bm{x}}(w)
=\{w\cdot g\mid g\in\langle\bm{x}\rangle\}.
\]
Thus $\mu_w(\bm{x})$ counts the entries of $\bm{x}$
belonging to the orbit of $w$ under the action
of $\langle\bm{x}\rangle$.
In this sense, the orbit multiplicities are
quandle analogues of the conjugacy-class
multiplicities used by Berger~\cite{Ber}.
\end{remark}

The following lemma characterizes trivial tuples in terms of their orbit multiplicities. 

\begin{lemma}\label{lem:triviality}
Let $\bm{x}\in X^m$ and $z\in X$.
Then $\bm{x}=(z,\ldots,z)$ if and only if
\[
\mu_w(\bm{x})=
\begin{cases}
m&\text{if }w=z,\\
0&\text{if }w\in X\setminus\{z\}.
\end{cases}
\]
In particular, if $\bm{y}\in X^m$ is trivial and
$\mu_w(\bm{y})=\mu_w(\bm{x})$ for every $w\in X$,
then $\bm{y}=\bm{x}$.
\end{lemma}

\begin{proof}
Suppose that $\bm{x}=(z,\ldots,z)$.
Since $z*z=z*^{-1}z=z$, we have $O_{\bm{x}}(z)=\{z\}$.

Suppose that $z\in O_{\bm{x}}(w)$ for some $w\in X$.
Then $w*^n z=z$ for some $n\in\mathbb{Z}$,
where $*^n z$ denotes the $n$-fold iteration
of the operation $*z$, with the inverse operation
used when $n<0$.
Since $z*z=z*^{-1}z=z$, we obtain
\[
w=(w*^n z)*^{-n}z=z*^{-n}z=z.
\]
Hence $z\notin O_{\bm{x}}(w)$ for every $w\ne z$.
Therefore
\[
\mu_w(\bm{x})=
\begin{cases}
m & \text{if }w=z,\\
0 & \text{if }w\ne z.
\end{cases}
\]

Conversely, suppose that the stated equalities hold.
Write $\bm{x}=(x_1,\ldots,x_m)$.
Since $x_i\in O_{\bm{x}}(x_i)$, we have
$\mu_{x_i}(\bm{x})>0$ for every $i$.
Hence $x_i=z$ for every $i$, so that
$\bm{x}=(z,\ldots,z)$.

The last assertion follows immediately
from the equivalence just proved.
\end{proof}

\begin{lemma}\label{lem:Berger-invariants}
If two tuples $\bm{x},\bm{y}\in X^m$ are Hurwitz equivalent,
then the following hold.
\begin{enumerate}
\item
$\Pi(\bm{x})=\Pi(\bm{y})$.
\item 
$O_{\bm{x}}(w)=O_{\bm{y}}(w)$ for every $w\in X$. 
\item
$\mu_w(\bm{x})=\mu_w(\bm{y})$ for every $w\in X$.
\end{enumerate}
\end{lemma}

\begin{proof}
Write $\bm{x}=(x_1,\ldots,x_m)$.
It suffices to consider the case
\[
\bm{y}=\bm{x}\cdot\sigma_i
=
(x_1,\ldots,x_{i-1},x_{i+1},x_i*x_{i+1},
x_{i+2},\ldots,x_m)
\]
for some $1\le i\le m-1$.

(i)
In $\As(X)$, we have
\[
x_{i+1}(x_i*x_{i+1})
=
x_{i+1}(x_{i+1}^{-1}x_i x_{i+1})
=
x_i x_{i+1}.
\]
Since all other factors are unchanged,
the total product is preserved.

(ii)
By the quandle axioms, for every $a\in X$ and
$\varepsilon\in\{1,-1\}$, we have 
\[
a*^\varepsilon x_i
=
((a*x_{i+1})*^\varepsilon(x_i*x_{i+1}))
*^{-1}x_{i+1} 
\]
and 
\[
a*^\varepsilon(x_i*x_{i+1})
=
((a*^{-1}x_{i+1})*^\varepsilon x_i)*x_{i+1}. 
\]
These identities show that the operations
determined by $x_i$ and $x_{i+1}$ generate
the same orbits as those determined by $x_i*x_{i+1}$ and $x_{i+1}$.
Hence, for every $w\in X$,
\[
O_{\bm{x}}(w)=O_{\bm{y}}(w).
\]

(iii) 
Let $w\in X$. 
By (ii), we have 
$x_i\in O_{\bm{x}}(w)$ if and only if $x_i*x_{i+1}\in O_{\bm{y}}(w)$. 
All other entries of $\bm{x}$ and $\bm{y}$ are the same. 
Therefore $\mu_w(\bm{x})=\mu_w(\bm{y})$. 
\end{proof}

In general, given $\bm{x},\bm{y}\in X^m$, 
it is not easy to determine whether
$\Pi(\bm{x})=\Pi(\bm{y})$. 
We therefore extract from $\Pi(\bm{x})$ a more tractable invariant as follows. 
For a group $G$ and a quandle homomorphism $\phi: X\to \Conj(G)$, 
define 
\[
\Pi_\phi(\bm{x})=\phi(x_1)\cdots\phi(x_m)\in G.
\]

\begin{lemma}
\label{lem:Pi_phi-invariant}
If two tuples $\bm{x},\bm{y}\in X^m$ are Hurwitz equivalent, 
then $\Pi_\phi(\bm{x})=\Pi_\phi(\bm{y})$ for every group $G$ and every quandle homomorphism $\phi:X\to\Conj(G)$.
\end{lemma}

\begin{proof}
By Lemma~\ref{lem:universality}, $\phi:X\to \Conj(G)$ induces a unique group homomorphism $\phi_\#: \As(X)\to G$ 
defined by $\phi_\#(x)=\phi(x)$ for $x\in X$. 
Using this homomorphism, we have 
$\Pi_\phi(\bm{x})=\phi_\#(\Pi(\bm{x}))$. 
By Lemma~\ref{lem:Berger-invariants}(i), we obtain 
\[
\Pi_\phi(\bm{x})=\phi_\#(\Pi(\bm{x}))=\phi_\#(\Pi(\bm{y}))=\Pi_\phi(\bm{y}). 
\qedhere 
\]
\end{proof}

\begin{proposition}\label{prop:Pi-vs-Pi_phi}
Let $\bm{x},\bm{y}\in X^m$. 
Then $\Pi(\bm{x})=\Pi(\bm{y})$ if and only if 
$\Pi_\phi(\bm{x})=\Pi_\phi(\bm{y})$ for every group $G$ and every quandle homomorphism $\phi:X\to\Conj(G)$. 
\end{proposition}

\begin{proof}
As in the proof of Lemma~\ref{lem:Pi_phi-invariant}, 
$\Pi(\bm{x})=\Pi(\bm{y})$ implies $\Pi_\phi(\bm{x})=\Pi_\phi(\bm{y})$, 
via the induced group homomorphism $\phi_\#$. 

Conversely, suppose that $\Pi_\phi(\bm{x})=\Pi_\phi(\bm{y})$ for every group $G$ and every quandle homomorphism $\phi:X\to\Conj(G)$. 
Let $G=\As(X)$ and let $\phi=\eta_X:X\to \Conj(\As(X))$ be the natural map. 
Then 
\[
\Pi(\bm{x})=\Pi_\phi(\bm{x})=\Pi_\phi(\bm{y})=\Pi(\bm{y}). 
\qedhere
\]
\end{proof}

Let $R_\infty=(\mathbb{Z},*)$ be the infinite dihedral quandle with 
\[
x*y=x*^{-1}y=2y-x
\quad (x,y\in\mathbb{Z}). 
\]
In~\cite{NNSW}, Nakamura, Nakanishi, Satoh, and the second author completely classified Hurwitz equivalence classes in $(R_\infty)^m$, 
using the following three invariants 
associated with $\bm{x}=(x_1,\ldots,x_m)\in(R_\infty)^m$. 

\begin{itemize}
\item 
$\Delta(\bm{x})=\displaystyle\sum_{i=1}^m(-1)^{i-1}x_i\in\mathbb{Z}$. 
\item 
$d(\bm{x})=\begin{cases}
\gcd\{x_i-x_1\mid 1\leq i\leq m\} & \text{if $\bm{x}$ is nontrivial}, \\
0 & \text{if $\bm{x}$ is trivial}. 
\end{cases}$

\item 
If $\bm{x}$ is nontrivial, 
let $M(\bm{x})$ be the multiset of residue classes $x_1,\ldots,x_m$ modulo $2d(\bm{x})$. 
If $\bm{x}$ is trivial, let $M(\bm{x})$ be the multiset of its entries. 
\end{itemize}

We now compare the three invariants
$\Delta(\bm{x})$, $d(\bm{x})$, and $M(\bm{x})$
with our two invariants $\Pi_\phi(\bm{x})$ and $\mu_w(\bm{x})$. 
Let $D_\infty=\langle r,s\mid s^2=(sr)^2=e\rangle$ 
be the infinite dihedral group. 
Define a quandle homomorphism
$\phi:R_\infty\to\Conj(D_\infty)$ by
\[
\phi(x)=sr^x
\quad (x\in R_\infty).
\]

\begin{theorem}\label{thm:Rinfty-invariants}
Let $\bm{x}=(x_1,\ldots,x_m)\in(R_\infty)^m$, and let
$\phi:R_\infty\to\Conj(D_\infty)$ be defined as above.
Then the following hold.
\begin{enumerate}
\item
We have
\[
\Pi_\phi(\bm{x})=
\begin{cases}
r^{-\Delta(\bm{x})}&\text{if $m$ is even},\\
sr^{\Delta(\bm{x})}&\text{if $m$ is odd}.
\end{cases}
\]
\item
If $d(\bm{x})>0$, then for every $1\le i\le m$, 
\[
O_{\bm{x}}(x_i)=x_i+2d(\bm{x})\Z, 
\]
and for every $w\in R_\infty$,
\[
\mu_w(\bm{x})
=\#\{i\mid x_i\equiv w\pmod{2d(\bm{x})}\}.
\]
\item
We have 
\[
\{w\in R_\infty\mid\mu_w(\bm{x})>0\}=
\begin{cases}
\{x_1\}&\text{if $d(\bm{x})=0$},\\
x_1+d(\bm{x})\Z&\text{if $d(\bm{x})>0$}.
\end{cases}
\]
\end{enumerate}
\end{theorem}

\begin{proof}
(i) 
Since $sr^ns=r^{-n}$ for every $n\in\mathbb{Z}$ in $D_\infty$, 
we have
\[
(sr^n)(sr^{n'})=r^{-n+n'}
\]
for every $n,n'\in\mathbb{Z}$. 

If $m$ is even, then
\begin{align*}
\Pi_\phi(\bm{x})
&=(sr^{x_1})(sr^{x_2})\cdots
(sr^{x_{m-1}})(sr^{x_m})\\
&=r^{-x_1+x_2-\cdots-x_{m-1}+x_m}
=r^{-\Delta(\bm{x})}.
\end{align*}
If $m$ is odd, then
\begin{align*}
\Pi_\phi(\bm{x})
&=(sr^{x_1})(sr^{x_2})(sr^{x_3})\cdots
(sr^{x_{m-1}})(sr^{x_m})\\
&=sr^{x_1-x_2+x_3-\cdots-x_{m-1}+x_m}
=sr^{\Delta(\bm{x})}.
\end{align*}

(ii)
Fix $1\le i\le m$.
Let $a\in O_{\bm{x}}(x_i)$. 
First, we show that $a-x_i\in 2d(\bm{x})\Z$
by induction on the minimum number of operations needed to obtain $a$ from $x_i$. 

If $a=x_i$, the assertion is clear.
Otherwise, we have $a=b*x_j$ or
$a=b*^{-1}x_j$ for some $1\leq j\leq m$ and 
some $b\in O_{\bm{x}}(x_i)$ obtained by fewer operations.
Since $b*^{-1}x_j=b*x_j$ in $R_\infty$, we have
\[
a-x_i=2x_j-b-x_i 
=2(x_j-x_i)-(b-x_i).
\]
By the induction hypothesis,
$b-x_i\in 2d(\bm{x})\Z$.
Since $d(\bm{x})$ divides $x_j-x_i$,
we obtain 
\[
a-x_i\in 2d(\bm{x})\Z.
\]
Hence $O_{\bm{x}}(x_i)
\subset x_i+2d(\bm{x})\Z$.

Conversely, let $a=x_i+2nd(\bm{x})$ for some $n\in\Z$.
By B\'ezout's identity, there exist integers
$c_2,\ldots,c_m$ such that
\[
d(\bm{x})=\sum_{j=2}^m c_j(x_j-x_1).
\]
Hence
\[
a=x_i+2n\sum_{j=2}^m c_j(x_j-x_1).
\]
For every $b\in R_\infty$ and $1\le j\le m$, we have
\[
(b*x_1)*x_j=b+2(x_j-x_1),
\qquad
(b*x_j)*x_1=b-2(x_j-x_1).
\]
By repeatedly applying these identities, 
using the appropriate one according to the sign of $nc_j$ for each $j$, 
we obtain $a$ from $x_i$. 
Hence $a\in O_{\bm{x}}(x_i)$. 
Combining this with the reverse inclusion, we obtain
\[
O_{\bm{x}}(x_i)=x_i+2d(\bm{x})\Z.
\]

For every $w\in R_\infty$, we have 
$x_i\in O_{\bm{x}}(w)$ if and only if 
$w\in O_{\bm{x}}(x_i)$. 
Therefore 
\[
\mu_w(\bm{x})
=\#\{i\mid x_i\equiv w
\pmod{2d(\bm{x})}\}.
\]

(iii)
If $d(\bm{x})=0$, then
$\bm{x}=(x_1,\ldots,x_1)$,
and the assertion follows from
Lemma~\ref{lem:triviality}.

Suppose that $d(\bm{x})>0$.
By (ii), we have
\[
\{w\in R_\infty\mid\mu_w(\bm{x})>0\}
=
\bigcup_{i=1}^m
\bigl(x_i+2d(\bm{x})\Z\bigr).
\]
Since $d(\bm{x})$ divides $x_i-x_1$ for every $i$,
we have
\[
\bigcup_{i=1}^m
\bigl(x_i+2d(\bm{x})\Z\bigr)
\subset x_1+d(\bm{x})\Z.
\]

On the other hand, 
let $a=x_1+nd(\bm{x})$ for some $n\in\Z$.
By the definition of $d(\bm{x})$,
the integers
\[
\frac{x_2-x_1}{d(\bm{x})},\ldots,
\frac{x_m-x_1}{d(\bm{x})}
\]
have greatest common divisor $1$.
Hence there exists $j$ such that
\[
x_j-x_1=(2p+1)d(\bm{x})
\]
for some $p\in\Z$.
If $n$ is even, then $a\in x_1+2d(\bm{x})\Z$. 
If $n$ is odd, then 
\[
a-x_j=(n-2p-1)d(\bm{x})
\in 2d(\bm{x})\Z,
\]
and hence $a\in x_j+2d(\bm{x})\Z$. 
Thus, in either case, we have
\[
a\in\bigcup_{i=1}^m
\bigl(x_i+2d(\bm{x})\Z\bigr).
\]
Therefore
\[
x_1+d(\bm{x})\Z
\subset
\bigcup_{i=1}^m
\bigl(x_i+2d(\bm{x})\Z\bigr).
\]
Combining the two inclusions gives the desired equality.
\end{proof}

As a consequence, we obtain the following.

\begin{corollary}\label{cor:Rinfty-invariants}
Let $\bm{x},\bm{y}\in(R_\infty)^m$, and let 
$\phi:R_\infty\to\Conj(D_\infty)$ be the quandle homomorphism defined above.
Then the following hold.
\begin{enumerate}
\item
$\Pi_\phi(\bm{x})=\Pi_\phi(\bm{y})$ if and only if
$\Delta(\bm{x})=\Delta(\bm{y})$.
\item
$\mu_w(\bm{x})=\mu_w(\bm{y})$ for every $w\in R_\infty$
if and only if
$d(\bm{x})=d(\bm{y})$ and $M(\bm{x})=M(\bm{y})$.
\end{enumerate}
\end{corollary}

\begin{proof}
(i) 
This follows from Theorem~\ref{thm:Rinfty-invariants}(i).

(ii) 
Write $\bm{x}=(x_1,\ldots,x_m)$ and $\bm{y}=(y_1,\ldots,y_m)$. 
Suppose that
$\mu_w(\bm{x})=\mu_w(\bm{y})$
for every $w\in R_\infty$.
Then
\[
\{w\in R_\infty\mid\mu_w(\bm{x})>0\}
=
\{w\in R_\infty\mid\mu_w(\bm{y})>0\}.
\]
By Theorem~\ref{thm:Rinfty-invariants}(iii),
we have $d(\bm{x})=d(\bm{y})$.

If $d(\bm{x})=d(\bm{y})=0$, then both tuples are trivial, and $x_1=y_1$.
Hence $M(\bm{x})=M(\bm{y})$. 
If $d(\bm{x})=d(\bm{y})>0$,
then Theorem~\ref{thm:Rinfty-invariants}(ii) gives
\[
\#\{i\mid x_i\equiv w\pmod{2d(\bm{x})}\}
=\#\{i\mid y_i\equiv w\pmod{2d(\bm{x})}\}
\]
for every $w\in R_\infty$.
Hence $M(\bm{x})=M(\bm{y})$.

Conversely, suppose that
$d(\bm{x})=d(\bm{y})$ and $M(\bm{x})=M(\bm{y})$.
If $d(\bm{x})=d(\bm{y})=0$,
then both tuples are trivial, and
$M(\bm{x})=M(\bm{y})$ implies $x_1=y_1$.
Hence $\mu_w(\bm{x})=\mu_w(\bm{y})$
for every $w\in R_\infty$.
If $d(\bm{x})=d(\bm{y})>0$,
then $M(\bm{x})=M(\bm{y})$ implies
\[
\#\{i\mid x_i\equiv w\pmod{2d(\bm{x})}\}
=
\#\{i\mid y_i\equiv w\pmod{2d(\bm{y})}\}
\]
for every $w\in R_\infty$.
By Theorem~\ref{thm:Rinfty-invariants}(ii),
we obtain
$\mu_w(\bm{x})=\mu_w(\bm{y})$
for every $w\in R_\infty$.
\end{proof}

\section{Trivial quandles and two sum constructions}
\label{sec:sum-constructions}

Let $X$ be a quandle.
For $\bm{x}=(x_1,\ldots,x_m)\in X^m$ and $w\in X$, set
\[
n_w(\bm{x})=\#\{i\mid x_i=w\}.
\]
Note that
\[
\sum_{w\in X}n_w(\bm{x})=m.
\]

We first classify Hurwitz equivalence classes for trivial quandles.
We then establish a classification theorem for a class of
amalgamated sums, which includes the quandles
$X_{3,2}$, $X_{4,2}$, $X_{4,3}$, and $X_{4,4}$.
Finally, we prove a reduction theorem for disjoint sums
of a quandle and a trivial quandle, showing that
Hurwitz equivalence for such a sum reduces to
Hurwitz equivalence for the former quandle
and the multiplicities of the entries in the latter.

\subsection{Trivial quandles}

Let $T_k=\{1,\ldots,k\}$ be the trivial quandle of order $k$, i.e.,
$x*y=x$ for all $x,y\in T_k$.

\begin{lemma}
\label{lem:trivial-normalform}
Every tuple $\bm{x}\in (T_k)^m$ is Hurwitz equivalent to 
\[
(\underbrace{1,\ldots,1}_{n_1(\bm{x})}, 
\underbrace{2,\ldots,2}_{n_2(\bm{x})}, \ldots, 
 \underbrace{k,\ldots,k}_{n_k(\bm{x})}).
\]
\end{lemma}

\begin{proof}
Write $\bm{x}=(x_1,\ldots,x_m)$. 
For every $j\in\{1,\ldots,m-1\}$, we have 
\[
(x_j,x_{j+1})\cdot \sigma_j=(x_{j+1},x_j).
\]
Thus the Hurwitz action factors through the natural action of the symmetric
group $\Sym_m$ by permutations of the entries. 
Hence every tuple is Hurwitz equivalent to the tuple displayed above.
\end{proof}

\begin{theorem}
\label{thm:trivial}
Let $\bm{x},\bm{y}\in (T_k)^m$. 
Then the following are equivalent.
\begin{enumerate}
\item $\bm{x}\sim\bm{y}$. 
\item 
$\mu_i(\bm{x})=\mu_i(\bm{y})$ for every $i\in\{1,\ldots,k\}$. 
\item $n_i(\bm{x})=n_i(\bm{y})$ for every $i\in\{1,\ldots,k\}$. 
\end{enumerate}
\end{theorem}

\begin{proof}
(i) $\Rightarrow$ (ii): 
This follows from Lemma~\ref{lem:Berger-invariants}(iii). 

(ii) $\Rightarrow$ (iii):
Since $T_k$ is a trivial quandle, we have
$i*x_j=i*^{-1}x_j=i$ for every $i\in T_k$ and every
entry $x_j$ of $\bm{x}$.
Hence $O_{\bm{x}}(i)=\{i\}$, 
and therefore 
$\mu_i(\bm{x})=n_i(\bm{x})$.
Similarly, $\mu_i(\bm{y})=n_i(\bm{y})$. 
Since $\mu_i(\bm{x})=\mu_i(\bm{y})$ by assumption, 
we obtain $n_i(\bm{x})=n_i(\bm{y})$ for every $i\in\{1,\ldots,k\}$.

(iii) $\Rightarrow$ (i): 
This follows from Lemma~\ref{lem:trivial-normalform}. 
\end{proof}

\subsection{A class of amalgamated sums}

We consider a special class of the amalgamated sums introduced by
Andruskiewitsch and Gra\~na~\cite{AG}.
Since their convention for the quandle operation is opposite to ours,
we describe the construction using the convention adopted in this paper.

Let $Y$ and $Z$ be disjoint nonempty finite sets,
regarded as trivial quandles.
For each $z\in Z$,
let $\theta_z\in\Sym(Y)$ be a permutation of $Y$,
and assume that
\[
\theta_z\theta_{z'}=\theta_{z'}\theta_z
\]
for all $z,z'\in Z$. 
Set $X=Y\sqcup Z$ and define a binary operation $*$ on $X$ by
\[
x*y=
\begin{cases}
x
& \text{if } y\in Y,\\
\theta_y(x)
& \text{if } x\in Y \text{ and } y\in Z,\\
x
& \text{if } x,y\in Z.
\end{cases}
\]
By \cite[Lemma~1.18]{AG},
this construction yields a quandle $X$,
which is an \emph{amalgamated sum} of $Y$ and $Z$.

\begin{example}\label{ex:small-amalgamated-sum}
The quandles $X_{3,2}$, $X_{4,2}$, $X_{4,3}$, and $X_{4,4}$
are amalgamated sums of the above form.
The corresponding sets $Y,Z$ and permutations $\theta_z$
are given as follows.
\begin{enumerate}
\item
For $X_{3,2}$, take
$Y=\{1,2\}$, $Z=\{3\}$, and $\theta_3=(1\,2)$.

\item
For $X_{4,2}$, take
$Y=\{1,2,4\}$, $Z=\{3\}$, and $\theta_3=(1\,2)$,
where $\theta_3$ fixes $4$.

\item
For $X_{4,3}$, take
$Y=\{1,2\}$, $Z=\{3,4\}$, and
$\theta_3=\theta_4=(1\,2)$.

\item
For $X_{4,4}$, take
$Y=\{1,2,3\}$, $Z=\{4\}$, and
$\theta_4=(1\,2\,3)$.
\end{enumerate}
These descriptions follow directly from the operation tables.
\end{example}

For $\bm{x}\in X^m$, let
\[
\Supp_Z(\bm{x})
=
\{z\in Z\mid n_z(\bm{x})>0\}
\]
and define $H(\bm{x})$ to be the subgroup of $\Sym(Y)$ 
generated by the permutations $\theta_z$ with $z\in\Supp_Z(\bm{x})$. 

\begin{lemma}\label{lem:amalgamated-sum-orbits}
Let $\bm{x}\in X^m$.
Then the following hold.
\begin{enumerate}
\item
For every $y\in Y$, 
$O_{\bm{x}}(y)
=\{h(y)\mid h\in H(\bm{x})\}$.

\item
For every $z\in Z$,
$O_{\bm{x}}(z)=\{z\}$.
\end{enumerate}
Consequently, for every $w\in X$, 
\[
\mu_w(\bm{x})=
\begin{cases}
\displaystyle\sum_{y'\in O_{\bm{x}}(w)}n_{y'}(\bm{x}) & \text{if }w\in Y,\\[1mm]
n_w(\bm{x}) & \text{if }w\in Z.
\end{cases}
\]
\end{lemma}

\begin{proof}
(i)
By the definition of the operation on $X$, 
for every $y\in Y$ and $x\in X$, we have 
\[
y*x=
\begin{cases}
y & \text{if }x\in Y,\\
\theta_x(y) & \text{if }x\in Z, 
\end{cases}
\qquad
y*^{-1}x=
\begin{cases}
y & \text{if }x\in Y,\\
\theta_x^{-1}(y) & \text{if }x\in Z. 
\end{cases}
\]
By the definition of $O_{\bm{x}}(y)$, its elements are obtained
by repeatedly applying the permutations $\theta_z$ and
$\theta_z^{-1}$ for $z\in\Supp_Z(\bm{x})$.
Hence
\[
O_{\bm{x}}(y)
=\{h(y)\mid h\in H(\bm{x})\}.
\]

(ii)
Since $z*x=z*^{-1}x=z$ for every $z\in Z$ and $x\in X$, 
every element of $Z$ is fixed under
the operations $*$ and $*^{-1}$ by the entries
of $\bm{x}$, and hence
\[
O_{\bm{x}}(z)=\{z\}.
\]

The formula for $\mu_w(\bm{x})$ follows immediately from (i), (ii), and the definition of $\mu_w$.
\end{proof}

Fix total orders on $Y$ and $Z$. 
Write $Y=\{y_1,\ldots,y_p\}$ and $Z=\{z_1,\ldots,z_q\}$ 
in increasing order. 

\begin{lemma}
\label{lem:amalgamated-sum-prenormalform}
Every tuple $\bm{x}\in X^m$ is Hurwitz equivalent to 
\[
(\underbrace{y_1,\ldots,y_1}_{n_{y_1}(\bm{x})},\ldots,
\underbrace{y_p,\ldots,y_p}_{n_{y_p}(\bm{x})}, 
\underbrace{z_1,\ldots,z_1}_{n_{z_1}(\bm{x})},\ldots,
\underbrace{z_q,\ldots,z_q}_{n_{z_q}(\bm{x})}).
\]
\end{lemma}

\begin{proof}
Since $Y$ and $Z$ are trivial quandles, 
we have 
\[
(y_i,y_j)\cdot\sigma_1=(y_j,y_i)
\]
for any $1\le i,j\le p$, and 
\[
(z_i,z_j)\cdot\sigma_1=(z_j,z_i)
\]
for any $1\le i,j\le q$. 
Moreover, since $x*y=x*^{-1}y=x$ for all $x\in X$ and $y\in Y$, 
we have 
\[
(y_i,z_j)\cdot\sigma_1^{-1}=(z_j,y_i),\qquad
(z_j,y_i)\cdot\sigma_1=(y_i,z_j)
\]
for any $1\le i\le p$ and $1\le j\le q$. 
Thus any two adjacent entries can be interchanged without changing their values. 
By repeatedly applying these interchanges, 
every tuple is Hurwitz equivalent to the tuple displayed above. 
\end{proof}

For $\bm{x}\in X^m$, write
\[
Y/H(\bm{x})=\{C_1,\ldots,C_r\},
\]
where the orbits are ordered according to their least elements,
and let $c_k$ denote the least element of $C_k$ $(k=1,\ldots,r)$.

\begin{proposition}\label{prop:amalgamated-sum-normalform}
Every tuple $\bm{x}\in X^m$ is Hurwitz equivalent to
\[
\bigl(
\underbrace{c_1,\ldots,c_1}_{\mu_{c_1}(\bm{x})},
\ldots,
\underbrace{c_r,\ldots,c_r}_{\mu_{c_r}(\bm{x})},
\underbrace{z_1,\ldots,z_1}_{\mu_{z_1}(\bm{x})},
\ldots,
\underbrace{z_q,\ldots,z_q}_{\mu_{z_q}(\bm{x})}
\bigr).
\]
\end{proposition}

\begin{proof}
By Lemma~\ref{lem:amalgamated-sum-prenormalform},
we have
\[
\bm{x}\sim
\bigl(
\underbrace{y_1,\ldots,y_1}_{n_{y_1}(\bm{x})},
\ldots,
\underbrace{y_p,\ldots,y_p}_{n_{y_p}(\bm{x})},
\underbrace{z_1,\ldots,z_1}_{n_{z_1}(\bm{x})},
\ldots,
\underbrace{z_q,\ldots,z_q}_{n_{z_q}(\bm{x})}
\bigr).
\]
If no entry of $\bm{x}$ belongs to $Y$, 
then the displayed tuple is already the desired normal form 
by Lemma~\ref{lem:amalgamated-sum-orbits}. 
If $\Supp_Z(\bm{x})=\varnothing$, then $H(\bm{x})=\{\id_Y\}$, 
and the same conclusion holds. 

Suppose that $n_{y_i}(\bm{x})>0$ and
$n_{z_j}(\bm{x})>0$
for some $1\le i\le p$ and $1\le j\le q$. 
As shown in the proof of Lemma~\ref{lem:amalgamated-sum-prenormalform}, 
we may move an entry equal to $z_j$ 
next to an entry equal to $y_i$ without changing any entries.
Then
\[
(y_i,z_j)\cdot\sigma_1^2
=
(\theta_{z_j}(y_i),z_j),
\qquad
(y_i,z_j)\cdot\sigma_1^{-2}
=
(\theta_{z_j}^{-1}(y_i),z_j).
\]
Thus an occurrence of $y_i$ can be replaced by
$\theta_{z_j}(y_i)$ or $\theta_{z_j}^{-1}(y_i)$
without changing the other entries.

Since these operations leave all entries in $Z$ unchanged,
the set $\Supp_Z(\bm{x})$ is preserved.
Moreover, $H(\bm{x})$ is generated by the permutations
$\theta_{z_j}$ with $n_{z_j}(\bm{x})>0$.
Hence, by repeatedly applying these operations,
we can replace each occurrence of $y_i$ by $c_k$,
where $C_k$ is the $H(\bm{x})$-orbit containing $y_i$.
Thus every entry belonging to $C_k$ can be replaced by $c_k$.

By Lemma~\ref{lem:amalgamated-sum-orbits},
the number of entries belonging to $C_k$ is
$\mu_{c_k}(\bm{x})$, while
\[
n_{z_j}(\bm{x})=\mu_{z_j}(\bm{x})
\qquad (1\le j\le q).
\]
Reordering the resulting entries gives the desired tuple.
\end{proof}

\begin{theorem}\label{thm:amalgamated-sum-classification}
Let $\bm{x},\bm{y}\in X^m$.
Then the following are equivalent.
\begin{enumerate}
\item
$\bm{x}\sim\bm{y}$.

\item
$\mu_w(\bm{x})=\mu_w(\bm{y})$ for every $w\in X$.
\end{enumerate}
\end{theorem}

\begin{proof}
(i) $\Rightarrow$ (ii):
This follows from Lemma~\ref{lem:Berger-invariants}(iii).

(ii) $\Rightarrow$ (i):
By Lemma~\ref{lem:amalgamated-sum-orbits},
we have
\[
n_{z_j}(\bm{x})=\mu_{z_j}(\bm{x})
=\mu_{z_j}(\bm{y})=n_{z_j}(\bm{y})
\]
for every $z_j\in Z$. 
Hence $\Supp_Z(\bm{x})=\Supp_Z(\bm{y})$,
which implies $H(\bm{x})=H(\bm{y})$. 
Thus $\bm{x}$ and $\bm{y}$ determine the same
orbit decomposition
\[
Y/H(\bm{x})=Y/H(\bm{y})
=\{C_1,\ldots,C_r\}.
\]

By Proposition~\ref{prop:amalgamated-sum-normalform},
$\bm{x}$ and $\bm{y}$ are Hurwitz equivalent
to normal forms defined using these orbits.
Since $\mu_w(\bm{x})=\mu_w(\bm{y})$ for every $w\in X$,
the corresponding blocks in these normal forms
have the same lengths.
Hence $\bm{x}$ and $\bm{y}$ are Hurwitz equivalent
to the same normal form, and therefore
$\bm{x}\sim\bm{y}$.
\end{proof}

We now apply the preceding classification theorem to four
nontrivial quandles of order at most four.

\begin{corollary}\label{cor:X32X42X43X44}
Let $X\in\{X_{3,2},X_{4,2},X_{4,3},X_{4,4}\}$, 
and let $\bm{x},\bm{y}\in X^m$.
Then the following are equivalent.
\begin{enumerate}
\item
$\bm{x}\sim\bm{y}$.

\item
$\mu_w(\bm{x})=\mu_w(\bm{y})$ for every $w\in X$.
\end{enumerate}
\end{corollary}

\begin{proof}
This follows from 
Example~\ref{ex:small-amalgamated-sum} and
Theorem~\ref{thm:amalgamated-sum-classification}.
\end{proof}

\subsection{Disjoint sums with trivial quandles}
\label{sec:disjoint-sum}

We consider a special case of the disjoint sum construction
of Andruskiewitsch and Gra\~na~\cite{AG},
using the convention for the quandle operation adopted
in this paper.

Let $Q$ be a quandle, and let
$Z$ be a nonempty finite set
disjoint from $Q$, regarded as a trivial quandle.
Set $X=Q\sqcup Z$ and define a binary operation $*$ on $X$ by
\[
x*y=
\begin{cases}
x*_Q y & \text{if }x,y\in Q,\\
x & \text{otherwise},
\end{cases}
\]
where $*_Q$ denotes the quandle operation of $Q$.
By \cite[Lemma~1.18]{AG},
this construction yields a quandle $X$,
which is the \emph{disjoint sum} of $Q$ and $Z$.
In particular, for every $z\in Z$ and $x\in X$,
we have
\[
x*z=x,
\qquad
z*x=z.
\]

For $\bm{x}\in X^m$, let $\widehat{\bm{x}}$
denote the tuple obtained from $\bm{x}$
by deleting all entries belonging to $Z$.
We allow $\widehat{\bm{x}}$ to be empty.
For tuples of length at most one,
Hurwitz equivalence means equality. 
For the empty tuple $\varnothing$, we set
\[
O_{\varnothing}(w)=\{w\},
\qquad
\mu_w(\varnothing)=0
\]
for every $w\in X$.

\begin{lemma}\label{lem:disjoint-sum-orbits}
Let $\bm{x}\in X^m$.
Then the following hold.
\begin{enumerate}
\item
For every $a\in Q$, $O_{\bm{x}}(a)
=O_{\widehat{\bm{x}}}(a)$. 
\item
For every $z\in Z$, 
$O_{\bm{x}}(z)=\{z\}$. 
\end{enumerate}
Consequently, for every $w\in X$,
\[
\mu_w(\bm{x})=
\begin{cases}
\mu_w(\widehat{\bm{x}}) & \text{if }w\in Q,\\
n_w(\bm{x}) & \text{if }w\in Z.
\end{cases}
\]
\end{lemma}

\begin{proof}
(i) 
By the definition of the operation on $X$,
for every $a\in Q$ and $x\in X$, we have
\[
a*x=
\begin{cases}
a*_Q x & \text{if }x\in Q,\\
a & \text{if }x\in Z,
\end{cases}
\qquad
a*^{-1}x=
\begin{cases}
a*^{-1}_Q x & \text{if }x\in Q,\\
a & \text{if }x\in Z.
\end{cases}
\]
Thus the entries belonging to $Z$ act trivially on $Q$,
and hence
\[
O_{\bm{x}}(a)
=
O_{\widehat{\bm{x}}}(a).
\]

(ii) 
Since $z*x=z*^{-1}x=z$ for every $z\in Z$, we have 
\[
O_{\bm{x}}(z)=\{z\}.
\]

The formula for $\mu_w(\bm{x})$ follows immediately 
from (i), (ii), and the definition of $\mu_w$.
\end{proof}

Fix a total order on $Z$, and 
write $Z=\{z_1,\ldots,z_q\}$ in increasing order. 

\begin{proposition}\label{prop:disjoint-sum-normalform}
Every tuple $\bm{x}\in X^m$ is Hurwitz equivalent to
\[
\bigl(
\widehat{\bm{x}},
\underbrace{z_1,\ldots,z_1}_{\mu_{z_1}(\bm{x})},
\ldots,
\underbrace{z_q,\ldots,z_q}_{\mu_{z_q}(\bm{x})}
\bigr).
\]
\end{proposition}

\begin{proof}
For every $x\in X$ and $z_j\in Z$, we have
\[
(z_j,x)\cdot\sigma_1=(x,z_j),
\qquad
(x,z_j)\cdot\sigma_1^{-1}=(z_j,x).
\]
Thus any two adjacent entries, at least one of which
belongs to $Z$, can be interchanged.
By repeatedly applying these interchanges,
we may move all entries belonging to $Z$ to the right
and arrange them in the order $z_1,\ldots,z_q$,
without changing the relative order of the entries in $Q$.
Since $\mu_{z_j}(\bm{x})=n_{z_j}(\bm{x})$ by Lemma~\ref{lem:disjoint-sum-orbits},
we obtain the desired tuple.
\end{proof}

\begin{theorem}\label{thm:disjoint-sum-classification}
Let $\bm{x},\bm{y}\in X^m$.
Then the following are equivalent.
\begin{enumerate}
\item
$\bm{x}\sim\bm{y}$.

\item
$\mu_{z_j}(\bm{x})=\mu_{z_j}(\bm{y})$ for every $z_j\in Z$,
and
$\widehat{\bm{x}}\sim\widehat{\bm{y}}$ in $Q$.
\end{enumerate}
\end{theorem}

\begin{proof}
(i) $\Rightarrow$ (ii):
By Lemma~\ref{lem:Berger-invariants}(iii),
we have
$\mu_{z_j}(\bm{x})=\mu_{z_j}(\bm{y})$
for every $z_j\in Z$.

To prove $\widehat{\bm{x}}\sim\widehat{\bm{y}}$ in $Q$,
it suffices to consider the action of 
$\sigma_i^{\pm1}$ on $\bm{x}=(x_1,\ldots,x_m)$.
If both $x_i$ and $x_{i+1}$ belong to $Q$,
let $k$ be the number of entries belonging to $Q$
among $x_1,\ldots,x_i$.
Then $x_i$ and $x_{i+1}$ are the $k$-th and
$(k+1)$-st entries of $\widehat{\bm{x}}$, respectively.
Thus the action of $\sigma_i^{\pm1}$ on $\bm{x}$
induces the action of $\sigma_k^{\pm1}$
on $\widehat{\bm{x}}$.
If at least one of $x_i$ and $x_{i+1}$
belongs to $Z$, then $\sigma_i^{\pm1}$
merely interchanges these entries,
so $\widehat{\bm{x}}$ remains unchanged.
Hence $\widehat{\bm{x}}\sim\widehat{\bm{y}}$ in $Q$.

(ii) $\Rightarrow$ (i):
By Proposition~\ref{prop:disjoint-sum-normalform},
$\bm{x}$ and $\bm{y}$ are Hurwitz equivalent to
\[
\bigl(
\widehat{\bm{x}},
\underbrace{z_1,\ldots,z_1}_{\mu_{z_1}(\bm{x})},
\ldots,
\underbrace{z_q,\ldots,z_q}_{\mu_{z_q}(\bm{x})}
\bigr), \qquad
\bigl(
\widehat{\bm{y}},
\underbrace{z_1,\ldots,z_1}_{\mu_{z_1}(\bm{y})},
\ldots,
\underbrace{z_q,\ldots,z_q}_{\mu_{z_q}(\bm{y})}
\bigr),
\]
respectively. 
Since $\mu_{z_j}(\bm{x})=\mu_{z_j}(\bm{y})$
for every $z_j\in Z$,
the corresponding terminal blocks are identical.
Moreover, $\widehat{\bm{x}}$ and
$\widehat{\bm{y}}$ have the same length
and are Hurwitz equivalent in $Q$ by assumption.
If their common length is at most one,
then $\widehat{\bm{x}}=\widehat{\bm{y}}$.
Otherwise, the braid realizing their
Hurwitz equivalence can be applied to
the initial blocks without changing
the terminal blocks.
Thus $\bm{x}\sim\bm{y}$.
\end{proof}

\section{The quandles $X_{3,3}$ and $X_{4,5}$}\label{sec:dihedral-order3}

\subsection{The dihedral quandle $X_{3,3}$ of order three}

Let 
\[
D_6=\langle r,s\mid r^3=s^2=(sr)^2=e\rangle
=\{e,r,r^2,s,sr,sr^2\} 
\]
be  the dihedral group of order six. 
Define a map $\phi:X_{3,3}\to\Conj(D_6)$ by 
\[
\phi(1)=s,\qquad 
\phi(2)=sr,\qquad
\phi(3)=sr^2. 
\]

\begin{lemma}\label{lem:X33-homomorphism}
The map $\phi:X_{3,3}\to\Conj(D_6)$ is a quandle homomorphism. 
\end{lemma}

\begin{proof}
Define a bijection $\psi:X_{3,3}\to R_3$ by 
\[
\psi(1)=0,\qquad
\psi(2)=1,\qquad
\psi(3)=2.
\]
A direct calculation shows that for $i,j\in X_{3,3}$,  
\[
\psi(i*j)=2\psi(j)-\psi(i).
\]
This shows that $\psi$ is an isomorphism between $X_{3,3}$ and $R_3$. 
Since the map $\phi$ is the composition of this isomorphism $\psi$ and 
the quandle homomorphism $\phi_3$ defined in Example~\ref{ex:dihedral-isomorphism}(i), 
$\phi$ is a quandle homomorphism. 
\end{proof}



\begin{lemma}\label{lem:X33-leading-one}
For any nontrivial tuple $\bm{x}\in(X_{3,3})^m$, 
there exists $\bm{x}'\in(X_{3,3})^{m-1}$ such that 
$\bm{x}\sim(1,\bm{x}')$. 
\end{lemma}

\begin{proof}
If $\bm{x}$ contains an entry equal to $1$, 
then we move this entry to the first position. 
Indeed, $(i,1)\cdot\sigma_1=(1,i*1)$ 
for every $i\in X_{3,3}$. 

Suppose that $\bm{x}$ has no entry equal to $1$. 
Since $\bm{x}$ is nontrivial, it contains both $2$ and $3$. 
Hence there is an adjacent pair of the form $(2,3)$ or $(3,2)$. 
Since 
\[(2,3)\cdot\sigma_1^{-1}=(1,2),\qquad 
(3,2)\cdot\sigma_1^{-1}=(1,3),\] 
the tuple $\bm{x}$ is Hurwitz equivalent to a nontrivial tuple containing an entry equal to $1$. 
We then move this entry equal to $1$ to the first position as above. 
\end{proof}

\begin{lemma}\label{lem:X33-case-trivial}
Let $z\in\{2,3\}$ and let $m\geq 3$. 
Then the following Hurwitz equivalence holds in $(X_{3,3})^m$:
\[
(1,\underbrace{z,\ldots,z}_{m-1})
\sim
\begin{cases}
(\underbrace{1,\ldots,1}_{m-1},z)
& \text{if $m$ is even}, \\
(\underbrace{1,\ldots,1}_{m-2},2,2)
& \text{if $m$ is odd}.
\end{cases}
\]
\end{lemma}

\begin{proof}
A direct calculation shows that for $z\in\{2,3\}$, 
\[
(1,z,z,z)\sim (1,1,1,z).
\]
Indeed, as illustrated in Figure~\ref{fig:pf-X33-case-trivial}, 
this equivalence is realized by 
\[
(1,z,z,z)\cdot
(\sigma_1\sigma_2^{-1}\sigma_1\sigma_3\sigma_2^{-1}\sigma_3)
=(1,1,1,z).\]  

\begin{figure}[htbp]
  \vspace{2mm}
  \begin{center}
   \hspace{1em}
    \begin{overpic}[]{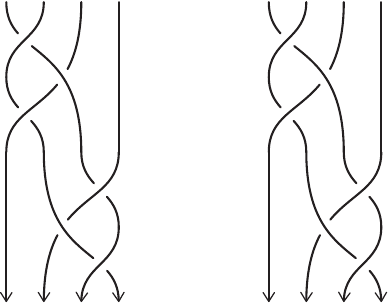}
      \put(0.5,150){$1$}
      \put(18.5,150){$2$}
      \put(36.5,150){$2$}
      \put(55,150){$2$}
      \put(23,121){$3$}
      \put(23,82){$3$}
      \put(0.5,-11){$1$}
      \put(18.5,-11){$1$}
      \put(36.5,-11){$1$}
      \put(55,-11){$2$}
      \put(127,150){$1$}
      \put(145,150){$3$}
      \put(163.25,150){$3$}
      \put(181,150){$3$}
      \put(150,121){$2$}
      \put(150,82){$2$}
      \put(127,-11){$1$}
      \put(145,-11){$1$}
      \put(163.25,-11){$1$}
      \put(181,-11){$3$}
    \end{overpic}
  \end{center}
  \vspace{2mm}
  \caption{The Hurwitz action of $\sigma_1\sigma_2^{-1}\sigma_1\sigma_3\sigma_2^{-1}\sigma_3$}
  \label{fig:pf-X33-case-trivial}
\end{figure}

Applying this equivalence repeatedly to consecutive blocks, we obtain
\[
(1,\underbrace{z,\ldots,z}_{m-1})
\sim
\begin{cases}
(\underbrace{1,\ldots,1}_{m-1},z)
& \text{if $m$ is even},\\
(\underbrace{1,\ldots,1}_{m-2},z,z)
& \text{if $m$ is odd}.
\end{cases}
\]

Thus the lemma is immediate when $m$ is even. 
When $m$ is odd, 
the case $z=2$ is already the desired form. 
For $z=3$, we use $(1,2,2)\sim (1,3,3)$, 
which follows from 
$(1,2,2)\cdot(\sigma_1\sigma_2^{-1}\sigma_1)=(1,3,3)$. 
Applying this equivalence to the last three entries gives
\[
(\underbrace{1,\ldots,1}_{m-2},3,3)
\sim
(\underbrace{1,\ldots,1}_{m-2},2,2).
\qedhere
\]
\end{proof}

\begin{proposition}\label{prop:X33-normalform}
Let $\bm{x}\in (X_{3,3})^m$ be nontrivial, 
and let $\phi:X_{3,3}\to\Conj(D_6)$ be the quandle homomorphism defined above. 
\begin{enumerate}
\item 
If $m=2$, then 
\[
\bm{x}\sim
\begin{cases}
(1,2) & \text{if } \Pi_\phi(\bm{x})=r,\\
(1,3) & \text{if } \Pi_\phi(\bm{x})=r^2.
\end{cases}
\]

\item
If $m\geq 3$ is odd, then
\[
\bm{x}\sim
\begin{cases}
(\underbrace{1,\ldots,1}_{m-1},2)
& \text{if } \Pi_\phi(\bm{x})=sr,\\
(\underbrace{1,\ldots,1}_{m-1},3)
& \text{if } \Pi_\phi(\bm{x})=sr^2,\\
(\underbrace{1,\ldots,1}_{m-2},2,2)
& \text{if } \Pi_\phi(\bm{x})=s.
\end{cases}
\]

\item
If $m\geq 4$ is even, then
\[
\bm{x}\sim
\begin{cases}
(\underbrace{1,\ldots,1}_{m-1},2)
& \text{if } \Pi_\phi(\bm{x})=r,\\
(\underbrace{1,\ldots,1}_{m-1},3)
& \text{if } \Pi_\phi(\bm{x})=r^2,\\
(\underbrace{1,\ldots,1}_{m-2},2,2)
& \text{if } \Pi_\phi(\bm{x})=e.
\end{cases}
\]
\end{enumerate}
\end{proposition}

\begin{proof}
(i) 
The nontrivial tuples of $(X_{3,3})^2$ are
\[
(1,2),\ (1,3),\ (2,1),\ (2,3),\ (3,1),\ (3,2).
\]
A direct calculation shows that  
\[
\Pi_\phi(1,2)=\Pi_\phi(2,3)=\Pi_\phi(3,1)=r, \qquad 
\Pi_\phi(1,3)=\Pi_\phi(2,1)=\Pi_\phi(3,2)=r^2.
\]
Moreover, since 
$(1,2)\cdot\sigma_1=(2,3)$ and 
$(2,3)\cdot\sigma_1=(3,1)$,
we have 
\[
(1,2)\sim (2,3)\sim (3,1).
\]
Similarly, we have $(1,3)\sim (3,2)\sim (2,1)$.

\medskip
(ii), (iii) 
We proceed by induction on $m$. 

First, suppose that $m=3$. 
By Lemma~\ref{lem:X33-leading-one}, 
we have $\bm{x}\sim (1,\bm{x}')$ for some $\bm{x}'\in (X_{3,3})^2$. 
Since $\bm{x}$ is nontrivial, $(1,\bm{x}')$ is also nontrivial. 

Since $\phi(1)=s$, we have
$\Pi_\phi(\bm{x})
=s\,\Pi_\phi(\bm{x}')$, 
and hence $\Pi_\phi(\bm{x}')=s\,\Pi_\phi(\bm{x})$. 
Thus
\[
\Pi_\phi(\bm{x}')=
\begin{cases}
r & \text{if } \Pi_\phi(\bm{x})=sr,\\
r^2 & \text{if } \Pi_\phi(\bm{x})=sr^2,\\
e & \text{if } \Pi_\phi(\bm{x})=s.
\end{cases}
\]

If $\Pi_\phi(\bm{x}')=r$ or $r^2$, then
assertion (i) implies that 
\[
\bm{x}'\sim
\begin{cases}
(1,2) & \text{if } \Pi_\phi(\bm{x})=sr,\\
(1,3) & \text{if } \Pi_\phi(\bm{x})=sr^2.
\end{cases}
\]
Therefore
\[
\bm{x}\sim
\begin{cases}
(1,1,2) & \text{if } \Pi_\phi(\bm{x})=sr,\\
(1,1,3) & \text{if } \Pi_\phi(\bm{x})=sr^2.
\end{cases}
\]

If $\Pi_\phi(\bm{x}')=e$, 
then $\bm{x}'$ is trivial. 
Write $\bm{x}'=(z,z)$ for some $z\in X_{3,3}$. 
Since $(1,\bm{x}')$ is nontrivial, we have $z\neq1$. 
Hence $z\in\{2,3\}$. 
Using $(1,2,2)\sim(1,3,3)$,
we obtain $\bm{x}\sim (1,2,2)$. 
This proves the case $m=3$.

Now, suppose that $m\geq 4$, 
and assume that the proposition holds for $m-1$. 
By Lemma~\ref{lem:X33-leading-one}, we have 
$\bm{x}\sim (1,\bm{x}')$ for some $\bm{x}'\in (X_{3,3})^{m-1}$. 
Again, we obtain 
$\Pi_\phi(\bm{x}')=s\,\Pi_\phi(\bm{x})$.

First, suppose that $\bm{x}'$ is nontrivial. 
If $m$ is even, then $m-1$ is odd. 
Hence, by the induction hypothesis applied to $\bm{x}'$, 
we have
\[
\bm{x}'\sim
\begin{cases}
(\underbrace{1,\ldots,1}_{m-2},2)
& \text{if } \Pi_\phi(\bm{x}')=sr,\\
(\underbrace{1,\ldots,1}_{m-2},3)
& \text{if } \Pi_\phi(\bm{x}')=sr^2,\\
(\underbrace{1,\ldots,1}_{m-3},2,2)
& \text{if } \Pi_\phi(\bm{x}')=s.
\end{cases}
\]
Since $\Pi_\phi(\bm{x}')=s\,\Pi_\phi(\bm{x})$, this gives
\[
\bm{x}\sim
\begin{cases}
(\underbrace{1,\ldots,1}_{m-1},2)
& \text{if } \Pi_\phi(\bm{x})=r,\\
(\underbrace{1,\ldots,1}_{m-1},3)
& \text{if } \Pi_\phi(\bm{x})=r^2,\\
(\underbrace{1,\ldots,1}_{m-2},2,2)
& \text{if } \Pi_\phi(\bm{x})=e.
\end{cases}
\]
This proves the even case when $\bm{x}'$ is nontrivial.

If $m$ is odd, then $m-1$ is even. 
Applying the induction hypothesis to $\bm{x}'$ gives, 
in the same way, 
\[
\bm{x}\sim
\begin{cases}
(\underbrace{1,\ldots,1}_{m-1},2)
& \text{if } \Pi_\phi(\bm{x})=sr,\\
(\underbrace{1,\ldots,1}_{m-1},3)
& \text{if } \Pi_\phi(\bm{x})=sr^2,\\
(\underbrace{1,\ldots,1}_{m-2},2,2)
& \text{if } \Pi_\phi(\bm{x})=s.
\end{cases}
\]
This proves the odd case when $\bm{x}'$ is nontrivial.

It remains to consider the case where $\bm{x}'$ is trivial. 
Then $\bm{x}'=(z,\ldots,z)$ for some $z\in X_{3,3}$. 
Since $(1,\bm{x}')$ is nontrivial, we have $z\neq1$. 
Hence 
\[
\bm{x}\sim
(1,\underbrace{z,\ldots,z}_{m-1})
\]
for some $z\in\{2,3\}$. 
If $m$ is even, then
\[
\Pi_\phi(1,\underbrace{2,\ldots,2}_{m-1})=r,
\qquad
\Pi_\phi(1,\underbrace{3,\ldots,3}_{m-1})=r^2.
\]
Since $\Pi_\phi$ is invariant under Hurwitz equivalence,
Lemma~\ref{lem:X33-case-trivial} implies 
\[
\bm{x}\sim 
\begin{cases}
(\underbrace{1,\ldots,1}_{m-1},2)
& \text{if } \Pi_\phi(\bm{x})=r,\\
(\underbrace{1,\ldots,1}_{m-1},3)
& \text{if } \Pi_\phi(\bm{x})=r^2.
\end{cases}
\]
If $m$ is odd, then
\[
\Pi_\phi(1,\underbrace{2,\ldots,2}_{m-1})
=
\Pi_\phi(1,\underbrace{3,\ldots,3}_{m-1})
=s.
\]
Thus $\Pi_\phi(\bm{x})=s$, and Lemma~\ref{lem:X33-case-trivial}
gives
\[
\bm{x}\sim(\underbrace{1,\ldots,1}_{m-2},2,2).
\]
This completes the induction.
\end{proof}

\begin{theorem}\label{thm:X33}
Let $\bm{x},\bm{y}\in(X_{3,3})^m$.
Then the following are equivalent.
\begin{enumerate}
\item $\bm{x}\sim\bm{y}$.
\item $\Pi_\phi(\bm{x})=\Pi_\phi(\bm{y})$ and
$\mu_i(\bm{x})=\mu_i(\bm{y})$ for every $i\in \{1,2,3\}$,
where $\phi:X_{3,3}\to\Conj(D_6)$ is defined above.
\end{enumerate}
\end{theorem}

\begin{proof}
(i) $\Rightarrow$ (ii):
This follows from Lemmas~\ref{lem:Berger-invariants}(iii) and~\ref{lem:Pi_phi-invariant}.

(ii) $\Rightarrow$ (i):
If $\bm{x}$ is trivial, Lemma~\ref{lem:triviality}
implies $\bm{x}=\bm{y}$.
If $\bm{x}$ is nontrivial, the same lemma implies
that $\bm{y}$ is nontrivial.
By Proposition~\ref{prop:X33-normalform},
$\bm{x}$ and $\bm{y}$ are then Hurwitz equivalent
to normal forms distinguished by $\Pi_\phi$.
Since $\Pi_\phi(\bm{x})=\Pi_\phi(\bm{y})$,
these normal forms coincide.
Hence $\bm{x}\sim\bm{y}$.
\end{proof}

\subsection{The quandle $X_{4,5}$}

Let $Q=\{1,2,3\}\subset X_{4,5}$.
By the operation table, $Q$ is a subquandle
isomorphic to $X_{3,3}$, and
\[
i*4=i,
\qquad
4*i=4
\qquad(i\in X_{4,5}).
\]
Thus $X_{4,5}$ is the disjoint sum of
$Q$ and the trivial quandle $Z=\{4\}$.

For $\bm{x}\in(X_{4,5})^m$, let
$\widehat{\bm{x}}$ denote the tuple obtained 
by deleting all entries equal to $4$. 
By Lemma~\ref{lem:disjoint-sum-orbits}, we have 
\[
\mu_i(\bm{x})=
\begin{cases}
\mu_i(\widehat{\bm{x}}) & \text{if }i\in Q,\\
n_4(\bm{x}) & \text{if }i=4.
\end{cases}
\]
Moreover, by Proposition~\ref{prop:disjoint-sum-normalform}, 
\[
\bm{x}\sim
\bigl(
\widehat{\bm{x}},
\underbrace{4,\ldots,4}_{\mu_4(\bm{x})}
\bigr). 
\]

Let $\phi:X_{3,3}\to\Conj(D_6)$ be the quandle homomorphism defined in the previous subsection.
Identifying $Q$ with $X_{3,3}$, define
$\phi':X_{4,5}\to\Conj(D_6)$ by 
\[
\phi'(i)=\phi(i)=sr^{i-1}\quad(i\in Q),
\qquad
\phi'(4)=e.
\]

\begin{lemma}\label{lem:X45-homomorphism}
The map $\phi'$ is a quandle homomorphism. 
Moreover, for every $\bm{x}\in(X_{4,5})^m$,
we have $\Pi_{\phi'}(\bm{x})
=\Pi_\phi(\widehat{\bm{x}})$. 
\end{lemma}

\begin{proof}
Since $\phi'|_Q=\phi$, 
Lemma~\ref{lem:X33-homomorphism} gives 
$\phi'(i*j)=\phi'(i)*\phi'(j)$ for every $i,j\in Q$. 

For every $i\in X_{4,5}$, we have $i*4=i$ and $4*i=4$. 
Hence  
\[
\phi'(i*4)=\phi'(i)=e\phi'(i)e=\phi'(i)*\phi'(4) 
\]
and
\[
\phi'(4*i)=\phi'(4)=\phi'(i)^{-1}e\phi'(i)=\phi'(4)*\phi'(i).
\]
Thus $\phi'$ is a quandle homomorphism.

Since $\phi'(4)=e$, deleting all entries equal to $4$ does not change the product $\Pi_{\phi'}(\bm{x})$.
Therefore $\Pi_{\phi'}(\bm{x})=\Pi_\phi(\widehat{\bm{x}})$. 
\end{proof}

\begin{theorem}\label{thm:X45}
Let $\bm{x},\bm{y}\in(X_{4,5})^m$.
Then the following are equivalent.
\begin{enumerate}
\item
$\bm{x}\sim\bm{y}$.

\item
$\Pi_{\phi'}(\bm{x})=\Pi_{\phi'}(\bm{y})$
and
$\mu_i(\bm{x})=\mu_i(\bm{y})$
for every $i\in\{1,2,3,4\}$, 
where $\phi':X_{4,5}\to \Conj(D_6)$ is the quandle homomorphism defined above. 
\end{enumerate}
\end{theorem}

\begin{proof}
(i) $\Rightarrow$ (ii):
This follows from Lemmas~\ref{lem:Berger-invariants}(iii) and~\ref{lem:Pi_phi-invariant}.

(ii) $\Rightarrow$ (i):
By Lemma~\ref{lem:disjoint-sum-orbits},
we have
\[
n_4(\bm{x})=\mu_4(\bm{x})
=\mu_4(\bm{y})=n_4(\bm{y}).
\]
Hence $\widehat{\bm{x}}$ and $\widehat{\bm{y}}$
have the same length.
Moreover, for every $i\in\{1,2,3\}$, we have 
\[
\mu_i(\widehat{\bm{x}})
=\mu_i(\bm{x})
=\mu_i(\bm{y})
=\mu_i(\widehat{\bm{y}}). 
\]
By Lemma~\ref{lem:X45-homomorphism}, we have 
\[
\Pi_\phi(\widehat{\bm{x}})
=\Pi_{\phi'}(\bm{x})
=\Pi_{\phi'}(\bm{y})
=\Pi_\phi(\widehat{\bm{y}}).
\]
If their common length is at most one, the equalities of
$\mu_i$ imply that $\widehat{\bm{x}}=\widehat{\bm{y}}$.
Otherwise, Theorem~\ref{thm:X33} gives
$\widehat{\bm{x}}\sim\widehat{\bm{y}}$.
By Theorem~\ref{thm:disjoint-sum-classification},
we conclude that $\bm{x}\sim\bm{y}$.
\end{proof}

\section{The dihedral quandle $X_{4,6}$ of order four}
\label{sec:dihedral-order4}

Let 
\[
D_8=\langle r,s\mid r^4=s^2=(sr)^2=e\rangle
=
\{e,r,r^2,r^3,s,sr,sr^2,sr^3\}
\]
be the dihedral group of order eight. 
Define a map
$\phi:X_{4,6}\to \Conj(D_8)$ by
\[
\phi(1)=s,\qquad
\phi(2)=sr^2,\qquad
\phi(3)=sr,\qquad
\phi(4)=sr^3.
\]

\begin{lemma}\label{lem:X46-homomorphism}
The map $\phi:X_{4,6}\to\Conj(D_8)$ is a quandle homomorphism.
\end{lemma}

\begin{proof}
Define a bijection $\psi:X_{4,6}\to R_4$ by
\[
\psi(1)=0,\qquad
\psi(2)=2,\qquad
\psi(3)=1,\qquad
\psi(4)=3.
\]
A direct calculation shows that for $i,j\in X_{4,6}$, 
\[
\psi(i*j)=2\psi(j)-\psi(i). 
\]
This shows that $\psi$ is an isomorphism between $X_{4,6}$ and $R_4$. 
Since $\phi=\phi_4\circ\psi$, the map $\phi$ is a quandle homomorphism. 
\end{proof}

Let $A=\{1,2\}$ and $B=\{3,4\}$. 
Then $A$ and $B$ are trivial subquandles of $X_{4,6}$.

\begin{lemma}\label{lem:X46-leading-one}
If $\bm{x}\in (X_{4,6})^m$ contains at least one entry of $A$ and 
at least one entry of $B$, 
then $\bm{x}\sim (1,\bm{x}')$ for some $\bm{x}'\in (X_{4,6})^{m-1}$ 
with one fewer entry from $A$ and the same number of entries from $B$.
\end{lemma}

\begin{proof}
If $\bm{x}$ contains an entry equal to $1$, 
then we move this entry to the first position. 
Indeed, $(i,1)\cdot\sigma_1=(1,i*1)$ for every $i\in X_{4,6}$, 
and $i*1$ belongs to $A$ if $i\in A$, and to $B$ if $i\in B$. 

Suppose that $\bm{x}$ has no entry equal to $1$. 
Then $\bm{x}$ contains an entry equal to $2$. 
Since $\bm{x}$ also contains an entry of $B$, 
there is an adjacent pair of the form $(2,b)$ or $(b,2)$ with $b\in B$. 
We have 
\[
(2,3)\cdot\sigma_1^2=(1,4), \quad 
(2,4)\cdot\sigma_1^2=(1,3), \quad
(3,2)\cdot\sigma_1^{-1}=(1,3), \quad 
(4,2)\cdot\sigma_1^{-1}=(1,4). 
\] 
Hence $\bm{x}$ is Hurwitz equivalent to a tuple 
containing an entry equal to $1$ and at least one entry of $B$. 
We then move this entry to the first position as above.
\end{proof}

\begin{lemma}\label{lem:X46-c-counts}
Let $\bm{x}\in (X_{4,6})^m$. 
If $\mu_1(\bm{x})+\mu_2(\bm{x})>0$ and $\mu_3(\bm{x})+\mu_4(\bm{x})>0$, 
then 
\[
\mu_1(\bm{x})=\mu_2(\bm{x})=n_1(\bm{x})+n_2(\bm{x}), \qquad 
\mu_3(\bm{x})=\mu_4(\bm{x})=n_3(\bm{x})+n_4(\bm{x}). 
\]
\end{lemma}

\begin{proof}
By the operation table of $X_{4,6}$ and the definition of the dual operation, 
for every $i\in X_{4,6}$, $a\in A$, and $b\in B$, 
we have 
\[a*i,\,a*^{-1}i\in A, \qquad 
b*i,\,b*^{-1}i\in B.\] 
Hence 
$O_{\bm{x}}(1),\,
O_{\bm{x}}(2)\subset A$ and 
$O_{\bm{x}}(3),\,
O_{\bm{x}}(4)\subset B$.

Since $\mu_1(\bm{x})+\mu_2(\bm{x})>0$, 
some entry of $\bm{x}$ belongs to $A$.
Similarly, since $\mu_3(\bm{x})+\mu_4(\bm{x})>0$, 
some entry of $\bm{x}$ belongs to $B$.
Thus $\bm{x}$ contains at least one entry
$a\in A$ and at least one entry $b\in B$.
By the operation table, $1*b=2$ and $2*b=1$ for every $b\in B$. 
Since $b$ occurs in $\bm{x}$, it follows that
\[
O_{\bm{x}}(1)
=O_{\bm{x}}(2)
=A.
\]
Similarly, since $a$ occurs in $\bm{x}$, we obtain 
\[
O_{\bm{x}}(3)
=O_{\bm{x}}(4)
=B.
\]
Therefore 
$\mu_1(\bm{x})=\mu_2(\bm{x})=n_1(\bm{x})+n_2(\bm{x})$ and 
$\mu_3(\bm{x})=\mu_4(\bm{x})=n_3(\bm{x})+n_4(\bm{x})$. 
\end{proof}

\begin{proposition}\label{prop:X46-normalform}
Every tuple $\bm{x}\in (X_{4,6})^m$  satisfies the following: 
\begin{enumerate}
\item 
If $\mu_3(\bm{x})=\mu_4(\bm{x})=0$, then 
\[
\bm{x}\sim
(\underbrace{1,\ldots,1}_{\mu_1(\bm{x})},\underbrace{2,\ldots,2}_{\mu_2(\bm{x})}).
\]

\item 
If $\mu_1(\bm{x})=\mu_2(\bm{x})=0$, then 
\[
\bm{x}\sim 
(\underbrace{3,\ldots,3}_{\mu_3(\bm{x})},\underbrace{4,\ldots,4}_{\mu_4(\bm{x})}). 
\]

\item
If $\mu_1(\bm{x})+\mu_2(\bm{x})>0$ and $\mu_3(\bm{x})+\mu_4(\bm{x})>0$, 
then $\bm{x}$ is Hurwitz equivalent to one of
\[
(\underbrace{1,\ldots,1}_{\mu_1(\bm{x})},\underbrace{3,\ldots,3}_{\mu_3(\bm{x})}), \qquad
(\underbrace{1,\ldots,1}_{\mu_1(\bm{x})},\underbrace{3,\ldots,3}_{\mu_3(\bm{x})-1},4).
\]
Moreover, the two normal forms are distinguished by the $\Pi_\phi$-invariant, 
where $\phi:X_{4,6}\to\Conj(D_8)$ is the quandle homomorphism defined above. 
\end{enumerate}
\end{proposition}

\begin{proof}
(i) 
Since $\mu_3(\bm{x})=\mu_4(\bm{x})=0$, we have
$n_3(\bm{x})=n_4(\bm{x})=0$.
Hence all entries of $\bm{x}$ lie in $A$.
Since $A$ is a trivial subquandle, the Hurwitz action merely
permutes the entries of $\bm{x}$.
Moreover, $O_{\bm{x}}(1)=\{1\}$ and $O_{\bm{x}}(2)=\{2\}$. 
Therefore $\mu_1(\bm{x})=n_1(\bm{x})$ and 
$\mu_2(\bm{x})=n_2(\bm{x})$. 
Thus, by permuting the entries, 
\[
\bm{x}\sim
(\underbrace{1,\ldots,1}_{\mu_1(\bm{x})},
 \underbrace{2,\ldots,2}_{\mu_2(\bm{x})}).
\]

(ii) 
The proof is identical to that of (i), 
with $A$, $1$, and $2$ replaced by $B$, $3$, and $4$, respectively. 

(iii) 
By Lemma~\ref{lem:X46-c-counts}, we may write 
\[
p=\mu_1(\bm{x})=n_1(\bm{x})+n_2(\bm{x})>0, \qquad
q=\mu_3(\bm{x})=n_3(\bm{x})+n_4(\bm{x})>0. 
\]
By repeatedly applying Lemma~\ref{lem:X46-leading-one}, 
we obtain 
\[
\bm{x}\sim
(\underbrace{1,\ldots,1}_{p},b_1,\ldots,b_q)
\] 
for some $b_1,\ldots,b_q\in B$. 
Since $B$ is a trivial subquandle, 
the entries $b_1,\ldots,b_q$ can be reordered arbitrarily. 
Thus
\[
\bm{x}\sim
(\underbrace{1,\ldots,1}_{p},
 \underbrace{3,\ldots,3}_{q-r},
 \underbrace{4,\ldots,4}_{r})
\]
for some integer $r$ with $0\leq r\leq q$. 

Since 
$(1,4,4)\cdot(\sigma_1\sigma_2^2\sigma_1)=(1,3,3)$, 
we have $(1,4,4)\sim(1,3,3)$. 
Since $p>0$, we can repeatedly apply this equivalence 
until at most one entry equal to $4$ remains. 
Therefore 
\[
\bm{x}\sim 
\begin{cases}
(\underbrace{1,\ldots,1}_{p},\underbrace{3,\ldots,3}_{q}) 
& \text{if $r$ is even}, \\
(\underbrace{1,\ldots,1}_{p},\underbrace{3,\ldots,3}_{q-1},4) 
& \text{if $r$ is odd}. 
\end{cases}
\]

Finally, these two normal forms are distinguished by the $\Pi_\phi$-invariant. 
Indeed, a direct calculation gives 
\[
\Pi_\phi(\underbrace{1,\ldots,1}_{p},\underbrace{3,\ldots,3}_{q})
=s^p(sr)^q, \qquad
\Pi_\phi(\underbrace{1,\ldots,1}_{p},\underbrace{3,\ldots,3}_{q-1},4)
=s^p(sr)^{q-1}(sr^3).
\]
If these two elements of $D_8$ were equal, 
then cancellation would give $sr=sr^3$, 
which is impossible in $D_8$. 
Hence
\[
\Pi_\phi(\underbrace{1,\ldots,1}_{p},\underbrace{3,\ldots,3}_{q})\neq \Pi_\phi(\underbrace{1,\ldots,1}_{p},\underbrace{3,\ldots,3}_{q-1},4).
\qedhere
\]
\end{proof}

\begin{theorem}\label{thm:X46}
Let $\bm{x},\bm{y}\in (X_{4,6})^m$. 
Then the following are equivalent. 

\begin{enumerate}
\item $\bm{x}\sim\bm{y}$.
\item 
$\Pi_\phi(\bm{x})=\Pi_\phi(\bm{y})$ and 
$\mu_i(\bm{x})=\mu_i(\bm{y})$ for every $i\in\{1,2,3,4\}$, 
where $\phi:X_{4,6}\to\Conj(D_8)$ is the quandle homomorphism defined above. 
\end{enumerate}
\end{theorem}

\begin{proof}
(i) $\Rightarrow$ (ii): 
This follows from Lemmas~\ref{lem:Berger-invariants}(iii) and~\ref{lem:Pi_phi-invariant}. 

(ii) $\Rightarrow$ (i): 
By Proposition~\ref{prop:X46-normalform}, 
$\bm{x}$ and $\bm{y}$ are Hurwitz equivalent to normal forms. 
These normal forms are distinguished by $\Pi_\phi$ and $(\mu_1,\mu_2,\mu_3,\mu_4)$. 
Since $\Pi_\phi(\bm{x})=\Pi_\phi(\bm{y})$ and $\mu_i(\bm{x})=\mu_i(\bm{y})$ for every $i\in\{1,2,3,4\}$, 
the two tuples are Hurwitz equivalent to the same normal form.
Hence $\bm{x}\sim\bm{y}$. 
\end{proof}

\section{The tetrahedral quandle $X_{4,7}$}\label{sec:tetrahedral}

Let 
\[
M=\mathbb F_2[t^{\pm1}]/(1+t+t^2)=\{0,1,t,1+t\} 
\]
be the quotient ring of the Laurent polynomial ring $\mathbb{F}_2[t^{\pm1}]$ by the ideal generated by $1+t+t^2$. 
Define a map 
$\omega:M\times M\to\mathbb{F}_2$ 
by
\[
\omega(a+bt,c+dt)=ac+bc+bd \quad(a,b,c,d\in\mathbb{F}_2).
\]

\begin{lemma}\label{lem:X47-property-omega}
The map $\omega$ satisfies the following: 
\begin{enumerate}
\item 
$\omega(x+x',y)=\omega(x,y)+\omega(x',y)$ and 
$\omega(x,y+y')=\omega(x,y)+\omega(x,y')$ for all $x,x',y,y'\in M$. 
\item 
$\omega(tx,ty)=\omega(x,y)=\omega(ty,x)$ for all $x,y\in M$. 
\item 
$\omega(t^{-1}x,t^{-1}x)=\omega(tx,x)$ for all $x\in M$. 
\end{enumerate}
\end{lemma}

\begin{proof}
Write 
\[x=a+bt,\qquad x'=a'+b't,\qquad 
y=c+dt,\qquad y'=c'+d't\] 
with $a,a',b,b',c,c',d,d'\in\mathbb F_2$. 
All computations below are in $\mathbb F_2$. 

(i) 
We have 
\begin{align*}
\omega(x+x',y)&=\omega((a+a')+(b+b')t,c+dt) \\
&=(a+a')c+(b+b')c+(b+b')d \\
&=(ac+bc+bd)+(a'c+b'c+b'd) \\
&=\omega(x,y)+\omega(x',y),  
\end{align*}
and similarly
\begin{align*}
\omega(x,y+y')&=\omega(a+bt,(c+c')+(d+d')t) \\
&=a(c+c')+b(c+c')+b(d+d') \\
&=(ac+bc+bd)+(ac'+bc'+bd') \\
&=\omega(x,y)+\omega(x,y'). 
\end{align*}

(ii) 
Since $1+t+t^2=0$ in $M$, we have $t^2=1+t$.
Hence 
\[
tx=t(a+bt)=b+(a+b)t, \qquad 
ty=t(c+dt)=d+(c+d)t. 
\]
Therefore  
\begin{align*}
\omega(tx,ty)
&=\omega\bigl(b+(a+b)t,\ d+(c+d)t\bigr)\\
&=bd+(a+b)d+(a+b)(c+d)\\
&=ac+bc+bd\\
&=\omega(x,y).
\end{align*}
Similarly, we obtain 
\begin{align*}
\omega(ty,x)
&=\omega\bigl(d+(c+d)t,\ a+bt\bigr)\\
&=da+(c+d)a+(c+d)b\\
&=ac+bc+bd\\
&=\omega(x,y). 
\end{align*}

(iii) 
Since $t^{-1}=1+t$, we have 
\[
t^{-1}x=(1+t)(a+bt)=(a+b)+at.
\]
Hence
\begin{align*}
\omega(t^{-1}x,t^{-1}x)
&=\omega\bigl((a+b)+at,\ (a+b)+at\bigr)\\
&=(a+b)^2+a(a+b)+a^2\\
&=a+b+ab.
\end{align*}
On the other hand, we have 
\begin{align*}
\omega(tx,x)
&=\omega\bigl(b+(a+b)t,\ a+bt\bigr)\\
&=ba+(a+b)a+(a+b)b\\
&=a+b+ab.
\end{align*}
Thus 
$\omega(t^{-1}x,t^{-1}x)=\omega(tx,x)$. 
\end{proof}

We now construct a group $G$ that will be used
to define a group-valued invariant $\Pi_\phi$ for $X_{4,7}$.
Its relation to the associated group of an
Alexander quandle, as described by
Clauwens~\cite{Cla}, will be explained below.
Set $G=\Z\times M\times\mathbb F_2$. 
For $(k,x,\alpha),(\ell,y,\beta)\in G$, define 
\[
(k,x,\alpha)(\ell,y,\beta)
=
(k+\ell,t^\ell x+y,\alpha+\beta+\omega(t^\ell x,y)).
\]

\begin{lemma}[{cf. \cite[Definition 1]{Cla}}]
\label{lem:X47-G-group}
With the multiplication above, $G$ is a group.
\end{lemma}

\begin{proof}
The identity element is $(0,0,0)$. 

For $(k,x,\alpha)\in G$, its inverse is
$(-k,t^{-k}x,\alpha+\omega(t^{-k}x,t^{-k}x))$. 
Indeed, we have 
\begin{align*}
(k,x,\alpha)(-k,t^{-k}x,\alpha+\omega(t^{-k}x,t^{-k}x)) 
&=(k-k, 2t^{-k}x, 2\alpha+2\omega(t^{-k}x,t^{-k}x)) \\
&=(0,0,0).
\end{align*}
The product in the reverse order is computed similarly. 

It remains to verify associativity. 
Let $(k,x,\alpha),(\ell,y,\beta),(m,z,\gamma)\in G$.
Associativity holds in the first two coordinates. 
For the third coordinate, the left-associated product gives 
\[
\alpha+\beta+\gamma
+\omega(t^\ell x,y)
+\omega(t^m(t^\ell x+y),z),
\]
whereas the right-associated product gives
\[
\alpha+\beta+\gamma
+\omega(t^m y,z)
+\omega(t^{\ell+m}x,t^m y+z).
\]
By Lemma~\ref{lem:X47-property-omega}(i), 
these are respectively equal to 
\[
\alpha+\beta+\gamma
+\omega(t^\ell x,y)
+\omega(t^{\ell+m}x,z)
+\omega(t^m y,z)
\]
and
\[
\alpha+\beta+\gamma
+\omega(t^m y,z)
+\omega(t^{\ell+m}x,t^m y)
+\omega(t^{\ell+m}x,z).
\]
By Lemma~\ref{lem:X47-property-omega}(ii), 
we have $\omega(t^\ell x,y)=\omega(t^{\ell+m}x,t^m y)$. 
Thus the two third coordinates are equal. 
Hence the multiplication is associative. 
\end{proof}

\begin{remark}\label{rem:Clauwens-group}
Let $M$ be an abelian group and $T$ an automorphism of $M$. 
Then $M$ admits a quandle structure defined by
\[
x*y=T(x)+y-T(y) \quad (x,y\in M).
\]
This quandle is called the \emph{Alexander quandle} 
associated with $M$ and $T$, and is denoted by $A(M,T)$. 
For example, the dihedral quandles and the tetrahedral quandle are Alexander quandles.
In~\cite{Cla}, Clauwens described the associated group of $A(M,T)$ in terms of $M$ and $T$. 
More precisely, define $\tau:M\otimes_{\Z}M\to M\otimes_{\Z}M$ by
\[
\tau(x\otimes y)=T(y)\otimes x,
\]
and set $S(M,T)=\operatorname{Coker}(1-\tau)$.
Clauwens proved that if $A(M,T)$ is connected, 
which is the case when $1-T$ is invertible, 
then $\As(A(M,T))$ is isomorphic to the group $\Z\times M\times S(M,T)$ with multiplication 
\[
(k,x,\alpha)(\ell,y,\beta)
=
(k+\ell,T^\ell(x)+y,
\alpha+\beta+[T^\ell(x)\otimes y]).
\]
In the case of the tetrahedral quandle, we have 
\[
M=\mathbb F_2[t^{\pm1}]/(1+t+t^2),
\]
where $T$ is multiplication by $t$.
Since $1-t=t^2$ in $M$, the map $1-T$ is invertible. 
By Lemma~\ref{lem:X47-property-omega}(i) and (ii), 
the map $\omega$ induces a homomorphism 
$\overline{\omega}:S(M,T)\to\mathbb F_2$ defined by 
\[
\overline{\omega}([x\otimes y])=\omega(x,y).
\]
The defining relations of $S(M,T)$ give 
\[
[1\otimes1]=[t\otimes1]=[t\otimes t],
\qquad
[1\otimes t]=0.
\]
Hence $S(M,T)$ is generated by $[1\otimes1]$. 
Since 
\[
\overline{\omega}([1\otimes1])=\omega(1,1)=1, 
\]
the homomorphism $\overline{\omega}$ is an isomorphism. 
Under this identification,
$[t^\ell x\otimes y]$ corresponds to
$\omega(t^\ell x,y)$.
Thus the group $G$ defined above is precisely
the group obtained from Clauwens' construction
in the case of the tetrahedral quandle.
\end{remark}

Define a bijection $\psi:X_{4,7}\to M$ by 
\[
\psi(1)=0,\qquad
\psi(2)=1,\qquad
\psi(3)=t,\qquad
\psi(4)=1+t. 
\]
A direct calculation shows that for $i,j\in X_{4,7}$,  
\[
\psi(i*j)=t\psi(i)+(1+t)\psi(j).
\]
This shows that $\psi$ is an isomorphism 
between $X_{4,7}$ and the tetrahedral quandle. 
Using this isomorphism, we 
define a map $\phi:X_{4,7}\to\Conj(G)$ by 
\[
\phi(i)=(1,\psi(i),0) \quad(i\in X_{4,7}). 
\]

\begin{lemma}\label{lem:X47-homomorphism}
The map $\phi:X_{4,7}\to\Conj(G)$ is a quandle homomorphism.
\end{lemma}

\begin{proof}
For $x,y\in M$, 
we first show that
\[
(1,y,0)^{-1}(1,x,0)(1,y,0)
=
(1,tx+(1+t)y,0).
\]
As in the proof of Lemma~\ref{lem:X47-G-group}, we have 
$(1,y,0)^{-1}
=(-1,t^{-1}y,\omega(t^{-1}y,t^{-1}y))$. 
Hence 
\begin{align*}
&(1,y,0)^{-1}(1,x,0)(1,y,0) \\
&=(0,x+y,\omega(t^{-1}y,t^{-1}y)+\omega(y,x))(1,y,0) \\
&=(1,t(x+y)+y, 
\omega(t^{-1}y,t^{-1}y)+\omega(y,x)+\omega(t(x+y),y)).
\end{align*}
The second coordinate is $t(x+y)+y=tx+(1+t)y$. 
The third coordinate is zero. 
Indeed, Lemma~\ref{lem:X47-property-omega} gives 
\begin{align*}
&\omega(t^{-1}y,t^{-1}y)+\omega(y,x)+\omega(t(x+y),y) \\
&=\omega(ty,y)+\omega(tx,y)+(\omega(tx,y)+\omega(ty,y))=0.
\end{align*}
Therefore $(1,y,0)^{-1}(1,x,0)(1,y,0)
=(1,tx+(1+t)y,0)$. 

It follows that, for $i,j\in X_{4,7}$,
\begin{align*}
\phi(i)*\phi(j)
&=\phi(j)^{-1}\phi(i)\phi(j) \\
&=(1,\psi(j),0)^{-1}(1,\psi(i),0)(1,\psi(j),0)\\
&=(1,t\psi(i)+(1+t)\psi(j),0)\\
&=(1,\psi(i*j),0)=\phi(i*j).
\qedhere
\end{align*}
\end{proof}

\begin{remark}
The proof of Theorem~1 in~\cite{Cla}
constructs an isomorphism
\[
\As(A(M,T))\to
\Z\times M\times S(M,T)
\]
that sends the generator corresponding to $x\in M$
to $(1,x,0)$.
Under the identification $S(M,T)\cong\mathbb F_2$ described in Remark~\ref{rem:Clauwens-group}, 
this gives an isomorphism
from $\As(A(M,T))$ to $G$.
Since
\[
(\phi\circ\psi^{-1})(x)=(1,x,0)
\qquad(x\in M),
\]
the induced group homomorphism
\[
(\phi\circ\psi^{-1})_\#:
\As(A(M,T))\to G
\]
coincides with this isomorphism.
Consequently, $\phi_\#:\As(X_{4,7})\to G$
is also an isomorphism.
\end{remark}

\begin{lemma}\label{lem:X47-leading-one}
Let $m\geq 3$. 
For any nontrivial tuple $\bm{x}\in (X_{4,7})^m$, 
there exists $\bm{x}'\in (X_{4,7})^{m-1}$ such that
$\bm{x}\sim(1,\bm{x}')$. 
\end{lemma}

\begin{proof}
If $\bm{x}$ contains an entry equal to $1$, 
then we move this entry to the first position. 
Indeed, $(i,1)\cdot\sigma_1=(1,i*1)$ for every $i\in X_{4,7}$. 

Suppose that $\bm{x}$ has no entry equal to $1$. 
Since $\bm{x}$ is nontrivial, 
it has a pair of adjacent distinct entries. 
Since $m\geq 3$, these two adjacent entries belong to a block of three consecutive entries. 
Therefore, $\bm{x}$ contains a nontrivial consecutive triple in $\{2,3,4\}^3$. 

We claim that every nontrivial triple in $\{2,3,4\}^3$ is Hurwitz equivalent to a triple containing $1$. 
Let $\rho:X_{4,7}\to X_{4,7}$ be a bijection defined by 
\[
\rho(1)=1,\qquad
\rho(2)=3,\qquad
\rho(3)=4,\qquad
\rho(4)=2.
\]
Since $\psi(\rho(i))=t\psi(i)$ for $i\in X_{4,7}$, 
and multiplication by $t$ is an automorphism of the tetrahedral quandle, 
the map $\rho$ is a quandle automorphism. 
We extend $\rho$ to $(X_{4,7})^3$ componentwise; 
that is, for $\bm{y}=(y_1,y_2,y_3)\in(X_{4,7})^3$, 
\[
\rho(\bm{y})=(\rho(y_1),\rho(y_2),\rho(y_3)).
\] 

There exist $24$ nontrivial triples in $\{2,3,4\}^3$. 
For the eight triples whose first entry is $2$, 
the following direct calculations hold:
\[
\begin{array}{@{}l@{\qquad}l@{}}
(2,2,3)\cdot(\sigma_2^2\sigma_1)=(4,1,2), &
(2,2,4)\cdot\sigma_2=(2,4,1),\\[2mm]
(2,3,2)\cdot\sigma_2=(2,2,1), &
(2,3,3)\cdot(\sigma_1\sigma_2)=(3,3,1),\\[2mm]
(2,3,4)\cdot(\sigma_2\sigma_1)=(4,1,2), &
(2,4,2)\cdot\sigma_1=(4,1,2),\\[2mm]
(2,4,3)\cdot\sigma_1=(4,1,3), &
(2,4,4)\cdot\sigma_1=(4,1,4).
\end{array}
\]
Since $\rho$ is a quandle automorphism, 
for every braid $\beta\in B_3$ and every
$\bm{y}\in (X_{4,7})^3$, we have 
$\rho(\bm{y}\cdot\beta)=\rho(\bm{y})\cdot\beta$. 
Therefore, applying $\rho$ and $\rho^2$ to both sides of 
each of the eight displayed equalities, 
we obtain analogous equivalences for all nontrivial triples in $\{2,3,4\}^3$ whose first entry is $3$ or $4$. 
Since $\rho$ fixes $1$, these equivalences also transform each such triple into a triple containing $1$. 
For example, the equality $(2,2,3)\cdot(\sigma_2^2\sigma_1)=(4,1,2)$ gives 
\[
(3,3,4)\cdot(\sigma_2^2\sigma_1)
=\rho(2,2,3)\cdot(\sigma_2^2\sigma_1)
=\rho(4,1,2)=(2,1,3)
\]
and 
\[
(4,4,2)\cdot(\sigma_2^2\sigma_1)
=\rho^2(2,2,3)\cdot(\sigma_2^2\sigma_1)
=\rho^2(4,1,2)=(3,1,4).
\]

Therefore, $\bm{x}$ is Hurwitz equivalent to a nontrivial tuple containing an entry equal to $1$. 
Moving this entry to the first position as above gives the desired form.
\end{proof}

For $m\geq2$, define the following tuples whenever all indicated block lengths are nonnegative:
\[
\begin{array}{@{}l@{\qquad}l@{}}
I_1(m)=(\underbrace{1,\ldots,1}_{m-1},2), & 
I_2(m)=(\underbrace{1,\ldots,1}_{m-2},2,2),\\
J_1(m)=(\underbrace{1,\ldots,1}_{m-1},3), & 
J_2(m)=(\underbrace{1,\ldots,1}_{m-2},3,3),\\
K_1(m)=(\underbrace{1,\ldots,1}_{m-1},4), & 
K_2(m)=(\underbrace{1,\ldots,1}_{m-2},4,4),\\
L_1(m)=(\underbrace{1,\ldots,1}_{m-2},2,3), & 
L_2(m)=(\underbrace{1,\ldots,1}_{m-3},2,2,2).
\end{array}
\]
Note that $I_2(2)$, $J_2(2)$, $K_2(2)$, and $L_2(3)$ are trivial tuples. 

\begin{lemma}\label{lem:X47-constant-tail}
For $m\geq3$, 
the following hold.
\[
(1,\underbrace{2,\ldots,2}_{m-1})
\sim
\begin{cases}
I_2(m) & \text{if } m\equiv0\pmod{3},\\
L_2(m) & \text{if } m\equiv1\pmod{3},\\
I_1(m) & \text{if } m\equiv2\pmod{3},
\end{cases}
\]
\[
(1,\underbrace{3,\ldots,3}_{m-1})
\sim
\begin{cases}
J_2(m) & \text{if } m\equiv0\pmod{3},\\
L_2(m) & \text{if } m\equiv1\pmod{3},\\
J_1(m) & \text{if } m\equiv2\pmod{3},
\end{cases}
\]
and
\[
(1,\underbrace{4,\ldots,4}_{m-1})
\sim
\begin{cases}
K_2(m) & \text{if } m\equiv0\pmod{3},\\
L_2(m) & \text{if } m\equiv1\pmod{3},\\
K_1(m) & \text{if } m\equiv2\pmod{3}.
\end{cases}
\]
\end{lemma}

\begin{proof}
A direct calculation shows that for $z\in\{2,3,4\}$, 
\[
(1,z,z,z,z)\sim(1,1,1,1,z). 
\]
Indeed, as illustrated in Figure~\ref{fig:pf-X47-constant-tail}, 
this equivalence is realized by 
\[
(1,z,z,z,z)\cdot
(\sigma_1\sigma_2\sigma_3\sigma_4^{-1}
 \sigma_3\sigma_2\sigma_1^{-1}\sigma_2\sigma_3\sigma_4)
=
(1,1,1,1,z). 
\] 

\begin{figure}[htbp]
  \vspace{2mm}
  \begin{center}
    \begin{overpic}[]{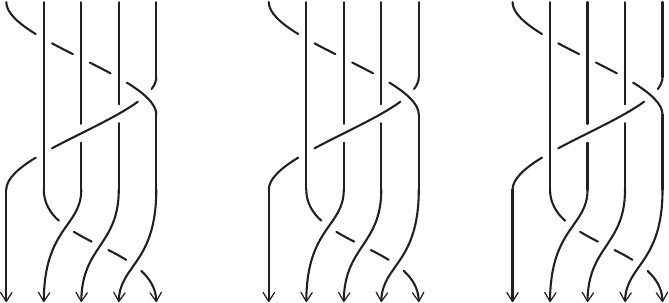}
      \put(0.75,149){$1$}
      \put(19,149){$2$}
      \put(37,149){$2$}
      \put(55,149){$2$}
      \put(73,149){$2$}
      \put(28,127){$4$}
      \put(46.5,118){$3$}
      \put(46.5,91){$4$}
      \put(40,36){$3$}
      \put(56,26.5){$4$}
      \put(0.75,-11){$1$}
      \put(19,-11){$1$}
      \put(37,-11){$1$}
      \put(55,-11){$1$}
      \put(73,-11){$2$}
      \put(127.25,149){$1$}
      \put(145.5,149){$3$}
      \put(163.5,149){$3$}
      \put(181.5,149){$3$}
      \put(199.5,149){$3$}
      \put(154,127){$2$}
      \put(172,118){$4$}
      \put(172,91){$2$}
      \put(167,36){$4$}
      \put(183,26.5){$2$}
      \put(127.25,-11){$1$}
      \put(145.5,-11){$1$}
      \put(163.5,-11){$1$}
      \put(181.5,-11){$1$}
      \put(199.5,-11){$3$}
      \put(244.5,149){$1$}
      \put(263,149){$4$}
      \put(281,149){$4$}
      \put(299,149){$4$}
      \put(317,149){$4$}
      \put(272,127){$3$}
      \put(289.5,118){$2$}
      \put(289.5,91){$3$}
      \put(284,36){$2$}
      \put(300,26.5){$3$}
      \put(244.5,-11){$1$}
      \put(263,-11){$1$}
      \put(281,-11){$1$}
      \put(299,-11){$1$}
      \put(317,-11){$4$}
    \end{overpic}
  \end{center}
  \vspace{2mm}
  \caption{The Hurwitz action of $\sigma_1\sigma_2\sigma_3\sigma_4^{-1}
 \sigma_3\sigma_2\sigma_1^{-1}\sigma_2\sigma_3\sigma_4$}
  \label{fig:pf-X47-constant-tail}
\end{figure}

We also have $(1,3,3,3)\sim(1,4,4,4)\sim(1,2,2,2)$, 
which follows from 
\[
(1,3,3,3)\cdot
(\sigma_1^{-1}\sigma_2\sigma_3^{-1}\sigma_2\sigma_1^{-1})
=
(1,2,2,2)
\]
and
\[
(1,4,4,4)\cdot
(\sigma_1\sigma_2\sigma_3^{-1}\sigma_2\sigma_1)
=
(1,2,2,2), 
\]
as illustrated in Figure~\ref{fig:pf-X47-constant-tail2}. 

\begin{figure}[htbp]
  \vspace{2mm}
  \begin{center}
    \begin{overpic}[]{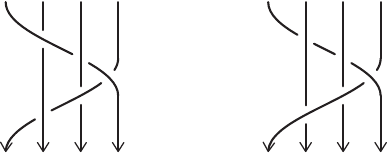}
      \put(0.5,77){$1$}
      \put(18.5,77){$3$}
      \put(36.5,77){$3$}
      \put(54.5,77){$3$}
      \put(45.5,20.5){$4$}
      \put(0.5,-11){$1$}
      \put(18.5,-11){$2$}
      \put(36.5,-11){$2$}
      \put(54.5,-11){$2$}
      \put(127,77){$1$}
      \put(144.5,77){$4$}
      \put(163,77){$4$}
      \put(181,77){$4$}
      \put(154.5,54.5){$3$}
      \put(127,-11){$1$}
      \put(144.5,-11){$2$}
      \put(163,-11){$2$}
      \put(181,-11){$2$}
    \end{overpic}
  \end{center}
  \vspace{2mm}
  \caption{The Hurwitz actions of $\sigma_1^{-1}\sigma_2\sigma_3^{-1}\sigma_2\sigma_1^{-1}$ and $\sigma_1\sigma_2\sigma_3^{-1}\sigma_2\sigma_1$}
  \label{fig:pf-X47-constant-tail2}
\end{figure}

Let $z\in\{2,3,4\}$. 
Repeatedly applying the first equivalence to a block consisting of one entry $1$ followed by four entries $z$, 
we obtain the following reductions:
\[
(1,\underbrace{z,\ldots,z}_{m-1})
\sim
\begin{cases}
(\underbrace{1,\ldots,1}_{m-2},z,z)
& \text{if } m\equiv0\pmod 3, \\
(\underbrace{1,\ldots,1}_{m-3},z,z,z)
& \text{if } m\equiv1\pmod 3,\\
(\underbrace{1,\ldots,1}_{m-1},z)
& \text{if } m\equiv2\pmod 3.
\end{cases}
\]

For $z=2$, this gives
\[
(1,\underbrace{2,\ldots,2}_{m-1})
\sim
\begin{cases}
I_2(m) & \text{if } m\equiv0\pmod 3, \\
L_2(m) & \text{if } m\equiv1\pmod 3, \\
I_1(m) & \text{if } m\equiv2\pmod 3.
\end{cases}
\]

For $z=3$, the cases $m\equiv0,2\pmod 3$ give $J_2(m)$ and $J_1(m)$,
respectively. If $m\equiv1\pmod 3$, then the above reduction gives
\[
(1,\underbrace{3,\ldots,3}_{m-1})
\sim
(\underbrace{1,\ldots,1}_{m-3},3,3,3).
\]
Using $(1,3,3,3)\sim(1,2,2,2)$ 
on the last four entries, we obtain
\[
(\underbrace{1,\ldots,1}_{m-3},3,3,3)
\sim
(\underbrace{1,\ldots,1}_{m-3},2,2,2)
=
L_2(m).
\]
Hence
\[
(1,\underbrace{3,\ldots,3}_{m-1})
\sim
\begin{cases}
J_2(m) & \text{if } m\equiv0\pmod 3,\\
L_2(m) & \text{if } m\equiv1\pmod 3,\\
J_1(m) & \text{if } m\equiv2\pmod 3.
\end{cases}
\]

The case $z=4$ is similar. 
The cases $m\equiv0,2\pmod 3$ give $K_2(m)$ and $K_1(m)$, respectively. 
If $m\equiv1\pmod 3$, then
\[
(1,\underbrace{4,\ldots,4}_{m-1})
\sim
(\underbrace{1,\ldots,1}_{m-3},4,4,4).
\]
Using $(1,4,4,4)\sim(1,2,2,2)$ 
on the last four entries, we obtain
\[
(\underbrace{1,\ldots,1}_{m-3},4,4,4)
\sim
(\underbrace{1,\ldots,1}_{m-3},2,2,2)
=
L_2(m).
\]
Therefore
\[
(1,\underbrace{4,\ldots,4}_{m-1})
\sim
\begin{cases}
K_2(m) & \text{if } m\equiv0\pmod 3,\\
L_2(m) & \text{if } m\equiv1\pmod 3,\\
K_1(m) & \text{if } m\equiv2\pmod 3.
\end{cases}
\qedhere
\]
\end{proof}

\begin{proposition}\label{prop:X47-normalform}
Every nontrivial tuple $\bm{x}\in (X_{4,7})^m$ satisfies the following:

\begin{enumerate}
\item
If $m=2$, then $\bm{x}$ is Hurwitz equivalent to one of
\[
I_1(2),\quad J_1(2),\quad K_1(2),\quad L_1(2).
\]

\item
If $m=3$, then $\bm{x}$ is Hurwitz equivalent to one of
\[
I_1(3),\ I_2(3),\ J_1(3),\ J_2(3),\
K_1(3),\ K_2(3),\ L_1(3).
\]

\item
If $m\geq4$, then $\bm{x}$ is Hurwitz equivalent to one of
\[
I_1(m),\ I_2(m),\ J_1(m),\ J_2(m),\
K_1(m),\ K_2(m),\ L_1(m),\ L_2(m).
\]
\end{enumerate}

Moreover, the listed normal forms are distinguished by the $\Pi_\phi$-invariant, 
where $\phi:X_{4,7}\to\Conj(G)$ is the quandle homomorphism defined above.
\end{proposition}

\begin{proof}
(i) 
The nontrivial pairs decompose into the following Hurwitz equivalence classes: 
\[
\begin{array}{@{}l@{\qquad}l@{}}
(1,2)\sim(2,4)\sim(4,1), &
(1,3)\sim(3,2)\sim(2,1), \\[2mm]
(1,4)\sim(4,3)\sim(3,1), &
(2,3)\sim(3,4)\sim(4,2).
\end{array}
\]
Thus every nontrivial pair is Hurwitz equivalent to one of
\[
I_1(2)=(1,2),\quad J_1(2)=(1,3),\quad K_1(2)=(1,4),\quad L_1(2)=(2,3).
\]

(ii) 
By Lemma~\ref{lem:X47-leading-one}, we have 
$\bm{x}\sim(1,\bm{x}')$ for some $\bm{x}'\in (X_{4,7})^2$. 
Since $\bm{x}$ is nontrivial, 
$(1,\bm{x}')$ is also nontrivial.

If $\bm{x}'$ is nontrivial, then the case $m=2$ implies that
$\bm{x}'$ is Hurwitz equivalent to one of
\[
I_1(2)=(1,2),\quad J_1(2)=(1,3),\quad K_1(2)=(1,4),\quad L_1(2)=(2,3).
\]
Hence $\bm{x}$ is Hurwitz equivalent to one of
\[
I_1(3)=(1,1,2),\quad J_1(3)=(1,1,3),\quad K_1(3)=(1,1,4),\quad L_1(3)=(1,2,3).
\]

If $\bm{x}'$ is trivial, then $\bm{x}'=(z,z)$ for some $z\in X_{4,7}$. 
Since $(1,\bm{x}')$ is nontrivial, we have $z\neq1$. 
Hence $z\in\{2,3,4\}$, and $\bm{x}$ is Hurwitz equivalent to one of 
\[
I_2(3)=(1,2,2),\quad
J_2(3)=(1,3,3),\quad
K_2(3)=(1,4,4).
\]

(iii) 
We proceed by induction on $m$, 
using the case $m=3$ proved above as the base case. 
By Lemma~\ref{lem:X47-leading-one}, we have 
$\bm{x}\sim(1,\bm{x}')$ for some $\bm{x}'\in (X_{4,7})^{m-1}$.

If $\bm{x}'$ is nontrivial, then the induction hypothesis shows that $\bm{x}'$ is Hurwitz equivalent to one of the listed normal forms of length $m-1$. 
Prepending $1$ to each of these normal forms gives the corresponding normal form of length $m$. 
Hence $(1,\bm{x}')$ is Hurwitz equivalent to one of the listed normal forms of length $m$. 

If $\bm{x}'$ is trivial, 
then $\bm{x}'=(z,\ldots,z)$ for some $z\in X_{4,7}$. 
Since $(1,\bm{x}')$ is nontrivial, 
we have $z\neq1$. 
Hence
\[
\bm{x}\sim(1,\underbrace{z,\ldots,z}_{m-1})
\]
for some $z\in\{2,3,4\}$. 
By Lemma~\ref{lem:X47-constant-tail}, this tuple is Hurwitz equivalent to one of the listed normal forms. 
This completes the induction. 

It remains to show that the listed normal forms are distinguished by the $\Pi_\phi$-invariant. 
A direct calculation using the multiplication in $G$ gives 
\[
\begin{array}{@{}l@{\qquad}l@{}}
\Pi_\phi(I_1(m))=(m,1,0), &
\Pi_\phi(I_2(m))=(m,1+t,1), \\[2mm]
\Pi_\phi(J_1(m))=(m,t,0), &
\Pi_\phi(J_2(m))=(m,1,1), \\[2mm]
\Pi_\phi(K_1(m))=(m,1+t,0), &
\Pi_\phi(K_2(m))=(m,t,1),\\[2mm]
\Pi_\phi(L_1(m))=(m,0,1), &
\Pi_\phi(L_2(m))=(m,0,0).
\end{array}
\]
Therefore, for each fixed $m$, the normal forms defined above have pairwise distinct values of the $\Pi_\phi$-invariant.
\end{proof}

\begin{theorem}\label{thm:X47}
Let $\bm{x},\bm{y}\in(X_{4,7})^m$.
Then the following are equivalent.
\begin{enumerate}
\item $\bm{x}\sim\bm{y}$.
\item $\Pi_\phi(\bm{x})=\Pi_\phi(\bm{y})$ and
$\mu_i(\bm{x})=\mu_i(\bm{y})$ for every $i\in X_{4,7}$,
where $\phi:X_{4,7}\to\Conj(G)$ is defined above.
\end{enumerate}
\end{theorem}

\begin{proof}
(i) $\Rightarrow$ (ii):
This follows from Lemmas~\ref{lem:Berger-invariants}(iii) and~\ref{lem:Pi_phi-invariant}.

(ii) $\Rightarrow$ (i):
If $\bm{x}$ is trivial, Lemma~\ref{lem:triviality}
implies $\bm{x}=\bm{y}$.
Otherwise, the same lemma shows that $\bm{y}$ is nontrivial.
By Proposition~\ref{prop:X47-normalform},
$\bm{x}$ and $\bm{y}$ are then Hurwitz equivalent
to normal forms distinguished by $\Pi_\phi$.
Their $\Pi_\phi$-invariants agree, so the normal forms coincide.
Hence $\bm{x}\sim\bm{y}$.
\end{proof}

\begin{proof}[Proof of Theorem~\ref{thm:main-classical}]
The only-if direction follows from Lemma~\ref{lem:Berger-invariants}(i) and (iii).

Conversely, suppose that
$\Pi(\bm{x})=\Pi(\bm{y})$ and
$\mu_w(\bm{x})=\mu_w(\bm{y})$ for every $w\in X$.
If $X$ is a trivial quandle or
$X\in\{X_{3,2},X_{4,2},X_{4,3},X_{4,4}\}$,
then Theorem~\ref{thm:trivial} and
Corollary~\ref{cor:X32X42X43X44} imply
$\bm{x}\sim\bm{y}$.

Now, suppose that $X\in\{X_{3,3},X_{4,5},X_{4,6},X_{4,7}\}$.
For the homomorphism $\phi$ (or~$\phi'$) used in the corresponding section,
Proposition~\ref{prop:Pi-vs-Pi_phi} gives
$\Pi_\phi(\bm{x})=\Pi_\phi(\bm{y})$ 
(or $\Pi_{\phi'}(\bm{x})=\Pi_{\phi'}(\bm{y})$).
Theorems~\ref{thm:X33}, \ref{thm:X45}, \ref{thm:X46},
and~\ref{thm:X47} then imply $\bm{x}\sim\bm{y}$.
\end{proof}

The complete invariants for Hurwitz equivalence
obtained in Sections~\ref{sec:sum-constructions}--\ref{sec:tetrahedral}
are summarized in Table~\ref{tab:classical-invariants}. 
In the last four cases, the homomorphisms $\phi$ and $\phi'$ are those defined in the corresponding sections.

\begin{table}[htbp]
\caption{Complete invariants for Hurwitz equivalence
in quandles of order at most four}
\centering
\small
\begin{tabular}{@{}ccc@{}}
\hline
Quandle $X$ & Complete invariants & Result \\
\hline
$X_{1,1},X_{2,1},X_{3,1},X_{4,1}$
& $(\mu_w)_{w\in X}$
& Theorem~\ref{thm:trivial} \\

$X_{3,2},X_{4,2},X_{4,3},X_{4,4}$
& $(\mu_w)_{w\in X}$
& Corollary~\ref{cor:X32X42X43X44} \\

$X_{3,3}$
& $\Pi_\phi,\,(\mu_w)_{w\in X}$
& Theorem~\ref{thm:X33} \\

$X_{4,5}$
& $\Pi_{\phi'},\,(\mu_w)_{w\in X}$
& Theorem~\ref{thm:X45} \\

$X_{4,6}$
& $\Pi_\phi,\,(\mu_w)_{w\in X}$
& Theorem~\ref{thm:X46} \\

$X_{4,7}$
& $\Pi_\phi,\,(\mu_w)_{w\in X}$
& Theorem~\ref{thm:X47} \\
\hline
\end{tabular}
\label{tab:classical-invariants}
\end{table}

\section{The Hurwitz action of the virtual braid group}
\label{sec:virtual}

Let $m\geq2$ be an integer. 
The \emph{virtual braid group $VB_m$} on $m$ strands is generated by 
$\sigma_1,\ldots,\sigma_{m-1}$ and $\tau_1,\ldots,\tau_{m-1}$
subject to the braid relations 
\[
\sigma_i\sigma_j=\sigma_j\sigma_i\quad (|i-j|>1),\qquad 
\sigma_i\sigma_{i+1}\sigma_i=\sigma_{i+1}\sigma_i\sigma_{i+1}, 
\] 
the permutation group relations 
\[
\tau_i^2=1,\qquad
\tau_i\tau_j=\tau_j\tau_i\quad (|i-j|>1), \qquad
\tau_i\tau_{i+1}\tau_i=\tau_{i+1}\tau_i\tau_{i+1}, 
\] 
and the mixed relations 
\[
\sigma_i\tau_j=\tau_j\sigma_i\quad (|i-j|>1), \qquad 
\tau_i\tau_{i+1}\sigma_i=\sigma_{i+1}\tau_i\tau_{i+1}.
\]
See~\cite{Kam07,KL} for more details. 

\begin{figure}[htbp]
  \vspace{2mm}
  \begin{center}
    \begin{overpic}[]{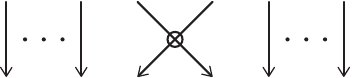}
      \put(-26,15.7){$\tau_{i}=$}
      \put(-2,43){$x_1$}
      \put(31.5,43){$x_{i-1}$}
      \put(60,43){$x_i$}
      \put(96,43){$x_{i+1}$}
      \put(122,43){$x_{i+2}$}
      \put(160,43){$x_{m}$}
      \put(-2,-11){$x_1$}
      \put(31.5,-11){$x_{i-1}$}
      \put(58,-11){$x_{i+1}$}
      \put(99,-11){$x_i$}
      \put(122,-11){$x_{i+2}$}
      \put(160,-11){$x_{m}$}
    \end{overpic}
  \end{center}
  \vspace{2mm}
  \caption{The Hurwitz action of $\tau_i$}
  \label{fig:Hurwitz-action-virtual}
\end{figure}

Let $X$ be a quandle. 
The \emph{Hurwitz action} of $VB_m$ on $X^m$ 
extends the Hurwitz action of $B_m$ defined in Section~\ref{sec:invariants} by specifying 
\[
(x_1,\ldots,x_i,x_{i+1},\ldots,x_m)\cdot \tau_i
 =
(x_1,\ldots,x_{i-1},x_{i+1},x_i,x_{i+2},\ldots,x_m), 
\]
see Figure~\ref{fig:Hurwitz-action-virtual}. 
Two tuples $\bm{x},\bm{y}\in X^m$ are \emph{virtually Hurwitz equivalent}
if there exists $\beta\in VB_m$ such that $\bm{x}\cdot\beta=\bm{y}$.
We write $\bm{x}\overset{v}{\sim}\bm{y}$ in this case. 
This equivalence relation is called \emph{virtual Hurwitz equivalence}.
Note that $\bm{x}\sim\bm{y}$ implies 
$\bm{x}\overset{v}{\sim}\bm{y}$.

\begin{lemma}\label{lem:virtual-invariants}
Let $X$ be a quandle and $\bm{x},\bm{y}\in X^m$.
If $\bm{x}\overset{v}{\sim}\bm{y}$, then
$\mu_w(\bm{x})=\mu_w(\bm{y})$ for every $w\in X$.
\end{lemma}

\begin{proof}
For the generators $\sigma_i$, this follows from Lemma~\ref{lem:Berger-invariants}(iii). 
The generator~$\tau_i$ merely interchanges two adjacent entries, 
so it leaves the operations defining $O_{\bm{x}}(w)$ unchanged, 
as well as the number of entries in each orbit.
Thus the orbit multiplicities are preserved by every generator.
\end{proof}

\begin{lemma}\label{lem:X33-virtual-normalform}
Every tuple $\bm{x}\in (X_{3,3})^m$ satisfies 
\[
\bm{x}\overset{v}{\sim}
\begin{cases}
(\underbrace{z,\ldots,z}_{m})
& \text{if }\bm{x}=(z,\ldots,z)\text{ is trivial},\\[1mm]
(\underbrace{1,\ldots,1}_{m-1},2)
& \text{if }\bm{x}\text{ is nontrivial}.
\end{cases}
\]
\end{lemma}

\begin{proof}
The assertion is immediate if $\bm{x}$ is trivial.
Suppose that $\bm{x}$ is nontrivial. 
By Proposition~\ref{prop:X33-normalform}, 
$\bm{x}$ is Hurwitz equivalent to one of 
\[
\begin{array}{ll}
(\underbrace{1,\ldots,1}_{m-1},2),\quad
(\underbrace{1,\ldots,1}_{m-1},3)
& (m\geq2),\\[2mm]
(\underbrace{1,\ldots,1}_{m-2},2,2)
& (m\geq3).
\end{array}
\]
For $m\geq2$, we have 
\[
(\underbrace{1,\ldots,1}_{m-1},3)\cdot(\tau_{m-1}\sigma_{m-1})
=
(\underbrace{1,\ldots,1}_{m-1},2),
\]
and, for $m\geq3$,
\[
(\underbrace{1,\ldots,1}_{m-2},2,2)\cdot(\tau_{m-2}\sigma_{m-2}\sigma_{m-1}\tau_{m-1})
=(\underbrace{1,\ldots,1}_{m-1},2).
\]
Thus $\bm{x}\overset{v}{\sim}(\underbrace{1,\ldots,1}_{m-1},2)$. 
\end{proof}

For $\bm{x}\in(X_{4,5})^m$, let
$\widehat{\bm{x}}$ denote the tuple obtained by deleting
all entries equal to~$4$.

\begin{lemma}\label{lem:X45-virtual-normalform}
Every tuple $\bm{x}\in(X_{4,5})^m$ satisfies
\[
\bm{x}\overset{v}{\sim}
\begin{cases}
(\underbrace{4,\ldots,4}_{m})
&\text{if }\widehat{\bm{x}}\text{ is empty},\\[1mm]
(\underbrace{z,\ldots,z}_{m-\mu_4(\bm{x})},
 \underbrace{4,\ldots,4}_{\mu_4(\bm{x})})
&\text{if }\widehat{\bm{x}}=(z,\ldots,z)\text{ is nonempty},\\[1mm]
(\underbrace{1,\ldots,1}_{m-\mu_4(\bm{x})-1},2,
 \underbrace{4,\ldots,4}_{\mu_4(\bm{x})})
&\text{if }\widehat{\bm{x}}\text{ is nontrivial}.
\end{cases}
\]
\end{lemma}

\begin{proof}
By Proposition~\ref{prop:disjoint-sum-normalform},
\[
\bm{x}\sim
(\widehat{\bm{x}},
 \underbrace{4,\ldots,4}_{\mu_4(\bm{x})}).
\]
If $\widehat{\bm{x}}$ is empty or trivial, this is
already the required form.
If $\widehat{\bm{x}}$ is nontrivial,
Lemma~\ref{lem:X33-virtual-normalform} gives
\[
\widehat{\bm{x}}\overset{v}{\sim}
(\underbrace{1,\ldots,1}_{m-\mu_4(\bm{x})-1},2).
\]
Applying the same virtual braid to the initial
block proves the assertion.
\end{proof}

\begin{lemma}\label{lem:X46-virtual-normalform}
Every tuple $\bm{x}\in(X_{4,6})^m$ satisfies 
\[
\bm{x}\overset{v}{\sim}
\begin{cases}
(\underbrace{1,\ldots,1}_{\mu_1(\bm{x})},
 \underbrace{2,\ldots,2}_{\mu_2(\bm{x})})
& \text{if } \mu_3(\bm{x})=\mu_4(\bm{x})=0,\\[1mm]
(\underbrace{3,\ldots,3}_{\mu_3(\bm{x})},
 \underbrace{4,\ldots,4}_{\mu_4(\bm{x})})
& \text{if } \mu_1(\bm{x})=\mu_2(\bm{x})=0,\\[1mm]
(\underbrace{1,\ldots,1}_{\mu_1(\bm{x})},
 \underbrace{3,\ldots,3}_{\mu_3(\bm{x})})
& \text{if } \mu_1(\bm{x})+\mu_2(\bm{x})>0
\text{ and } \mu_3(\bm{x})+\mu_4(\bm{x})>0.
\end{cases}
\]
\end{lemma}

\begin{proof}
The first two cases follow from Proposition~\ref{prop:X46-normalform}. 
Suppose that
$\mu_1(\bm{x})+\mu_2(\bm{x})>0$ and 
$\mu_3(\bm{x})+\mu_4(\bm{x})>0$.
By Lemma~\ref{lem:X46-c-counts}, we have 
\[
\mu_1(\bm{x})=\mu_2(\bm{x})>0,
\qquad
\mu_3(\bm{x})=\mu_4(\bm{x})>0.
\]
Proposition~\ref{prop:X46-normalform} implies that $\bm{x}$ is Hurwitz
equivalent to either
\[
(\underbrace{1,\ldots,1}_{\mu_1(\bm{x})},
 \underbrace{3,\ldots,3}_{\mu_3(\bm{x})})
\quad \text{or}\quad
(\underbrace{1,\ldots,1}_{\mu_1(\bm{x})},
 \underbrace{3,\ldots,3}_{\mu_3(\bm{x})-1},4).
\]
In the latter case, 
since $(1,4)\cdot(\tau_1\sigma_1)=(1,3)$ and 
$B=\{3,4\}$ is a trivial subquandle of $X_{4,6}$, 
we have
\[
(\underbrace{1,\ldots,1}_{\mu_1(\bm{x})}, 
\underbrace{3,\ldots,3}_{\mu_3(\bm{x})-1},4)
\sim 
(\underbrace{1,\ldots,1}_{\mu_1(\bm{x})},4, 
\underbrace{3,\ldots,3}_{\mu_3(\bm{x})-1})
\overset{v}{\sim}
(\underbrace{1,\ldots,1}_{\mu_1(\bm{x})}, 
\underbrace{3,\ldots,3}_{\mu_3(\bm{x})}).
\qedhere
\]
\end{proof}

\begin{lemma}\label{lem:X47-virtual-normalform}
Every tuple $\bm{x}\in(X_{4,7})^m$ satisfies 
\[
\bm{x}\overset{v}{\sim}
\begin{cases}
(\underbrace{z,\ldots,z}_{m})
& \text{if }\bm{x}=(z,\ldots,z)\text{ is trivial},\\[1mm]
(\underbrace{1,\ldots,1}_{m-1},2)
& \text{if }\bm{x}\text{ is nontrivial}.
\end{cases}
\]
\end{lemma}

\begin{proof}
The assertion is immediate if $\bm{x}$ is trivial.
Suppose that $\bm{x}$ is nontrivial. 
By Proposition~\ref{prop:X47-normalform}, 
$\bm{x}$ is Hurwitz equivalent
to one of
\[
I_1(m),\ I_2(m),\ J_1(m),\ J_2(m),
\ K_1(m),\ K_2(m),\ L_1(m),\ L_2(m),
\]
whenever these normal forms are defined.

We show that all these normal forms are virtually Hurwitz equivalent to
\[
(\underbrace{1,\ldots,1}_{m-1},2).
\]
First, the operation table of $X_{4,7}$ gives
\[
(1,3)\cdot(\tau_1\sigma_1)=(1,4),
\qquad
(1,4)\cdot(\tau_1\sigma_1)=(1,2), \qquad
(2,3)\cdot(\tau_1\sigma_1\tau_1)=(1,2).
\]
Therefore
\[
J_1(m)\overset{v}{\sim} K_1(m)\overset{v}{\sim} I_1(m),
\qquad
L_1(m)\overset{v}{\sim} I_1(m).
\]

Moreover, as illustrated in Figure~\ref{fig:pf-X47-virtual-normalform}, 
we have 
\[
(1,z,z)\cdot(\tau_1\sigma_1\sigma_2\tau_2)=(1,1,z) 
\] 
for $z\in\{2,3,4\}$. Hence
\[
I_2(m)\overset{v}{\sim} I_1(m), 
\qquad
J_2(m)\overset{v}{\sim} J_1(m)\overset{v}{\sim} I_1(m), 
\qquad
K_2(m)\overset{v}{\sim} K_1(m)\overset{v}{\sim} I_1(m).
\]

\begin{figure}[htbp]
  \vspace{2mm}
  \begin{center}
    \begin{overpic}[]{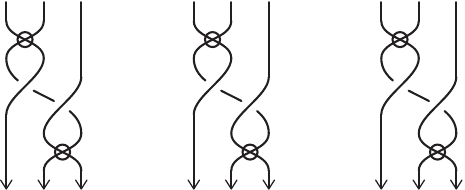}
      \put(0.5,96){$1$}
      \put(19.5,96){$2$}
      \put(36.5,96){$2$}
      \put(21,48){$3$}
      \put(0.5,-11){$1$}
      \put(19.5,-11){$1$}
      \put(36.5,-11){$2$}
      \put(91,96){$1$}
      \put(109,96){$3$}
      \put(127,96){$3$}
      \put(111,48){$4$}
      \put(91,-11){$1$}
      \put(109,-11){$1$}
      \put(127,-11){$3$}
      \put(181.5,96){$1$}
      \put(199.5,96){$4$}
      \put(217,96){$4$}
      \put(202,48){$2$}
      \put(181.5,-11){$1$}
      \put(199.5,-11){$1$}
      \put(217,-11){$4$}
    \end{overpic}
  \end{center}
  \vspace{2mm}
  \caption{The Hurwitz action of $\tau_1\sigma_1\sigma_2\tau_2$}
  \label{fig:pf-X47-virtual-normalform}
\end{figure}

Finally, if $\bm{x}$ is Hurwitz equivalent to 
\[
L_2(m)=
(\underbrace{1,\ldots,1}_{m-3},2,2,2), 
\]
then since $\bm{x}$ is nontrivial, we have $m\geq4$. 
Applying the equivalence
\[
(1,2,2)\overset{v}{\sim}(1,1,2)
\]
twice to the terminal part, we obtain
\[
L_2(m)\overset{v}{\sim}
(\underbrace{1,\ldots,1}_{m-1},2)
=
I_1(m).
\]
Thus every nontrivial tuple $\bm{x}\in(X_{4,7})^m$ is 
virtually Hurwitz equivalent to
\[
(\underbrace{1,\ldots,1}_{m-1},2).
\qedhere
\]
\end{proof}

\begin{theorem}\label{thm:main-virtual-X33X45X46X47}
Let $X\in\{X_{3,3},X_{4,5},X_{4,6},X_{4,7}\}$,
and let $\bm{x},\bm{y}\in X^m$.
Then $\bm{x}\overset{v}{\sim}\bm{y}$ if and only if
$\mu_w(\bm{x})=\mu_w(\bm{y})$ for every $w\in X$.
\end{theorem}

\begin{proof}
The only-if direction follows from
Lemma~\ref{lem:virtual-invariants}.
Conversely, suppose that $\mu_w(\bm{x})=\mu_w(\bm{y})$ for every $w\in X$.

First, suppose that $X=X_{3,3}$ or $X=X_{4,7}$.
If $\bm{x}$ is trivial, Lemma~\ref{lem:triviality}
gives $\bm{x}=\bm{y}$.
Otherwise, the same lemma implies that $\bm{y}$ is nontrivial.
The respective Lemmas~\ref{lem:X33-virtual-normalform}
and~\ref{lem:X47-virtual-normalform} then show that
both tuples are virtually Hurwitz equivalent to the same normal form.

If $X=X_{4,5}$, Lemma~\ref{lem:disjoint-sum-orbits} gives $\mu_4(\bm{x})=n_4(\bm{x})$ and $\mu_i(\widehat{\bm{x}})=\mu_i(\bm{x})$ for $i\in\{1,2,3\}$, 
and similarly for $\bm{y}$. 
These equalities determine whether $\widehat{\bm{x}}$ and $\widehat{\bm{y}}$
are empty, trivial, or nontrivial. 
In the trivial case, they also determine their common entry. 
Hence Lemma~\ref{lem:X45-virtual-normalform} shows that $\bm{x}$ and $\bm{y}$ are virtually Hurwitz equivalent to the same normal form. 

Finally, if $X=X_{4,6}$, the normal forms in Lemma~\ref{lem:X46-virtual-normalform} are determined by the values $\mu_i$. 
Hence $\bm{x}\overset{v}{\sim}\bm{y}$. 
\end{proof}

\begin{proof}[Proof of Theorem~\ref{thm:main-virtual}]
The only-if direction follows from
Lemma~\ref{lem:virtual-invariants}.
Conversely, suppose that
$\mu_w(\bm{x})=\mu_w(\bm{y})$ for every $w\in X$.
If $X$ is a trivial quandle or
$X\in\{X_{3,2},X_{4,2},X_{4,3},X_{4,4}\}$,
Theorem~\ref{thm:trivial} and
Corollary~\ref{cor:X32X42X43X44} give
$\bm{x}\sim\bm{y}$, hence
$\bm{x}\overset{v}{\sim}\bm{y}$.
For the remaining quandles, the result follows from
Theorem~\ref{thm:main-virtual-X33X45X46X47}.
\end{proof}

\begin{remark}
The \emph{welded braid group} $WB_m$~\cite{Kam07,KL}, 
also called the \emph{braid-permutation group}~\cite{FRR}, 
is the quotient of the virtual braid group $VB_m$ by the additional relations 
\[
\tau_i\sigma_{i+1}\sigma_i
=
\sigma_{i+1}\sigma_i\tau_{i+1}
\quad (1\leq i\leq m-2).
\]
For a quandle $X$, the Hurwitz action of $VB_m$ on $X^m$ factors through this quotient. 
Indeed, for any $x,y,z\in X$, we have
\[
(x,y,z)\cdot(\tau_1\sigma_2\sigma_1)
=
(z,y*z,x*z)
\]
and
\[
(x,y,z)\cdot(\sigma_2\sigma_1\tau_2)
=
(z,y*z,x*z).
\]
Applying this calculation to the entries in position $i$, $i+1$, and $i+2$ 
shows that the additional defining relations of $WB_m$ act trivially on
$X^m$. 
Hence the Hurwitz action of $VB_m$ on $X^m$ induces a well-defined Hurwitz action of $WB_m$ on $X^m$. 
Consequently, the equivalence relation on $X^m$ induced by the action of
$WB_m$ coincides with virtual Hurwitz equivalence.
\end{remark}



\begin{thebibliography}{99}
\bibitem{AG}
N. Andruskiewitsch and M. Gra\~na, 
\textit{From racks to pointed Hopf algebras}, 
Adv. Math. \textbf{178} (2003), no.~2, 177--243. 

\bibitem{BT}
T. Ben-Itzhak and M. Teicher, 
\textit{Graph theoretic method for determining Hurwitz equivalence 
in the symmetric group}, 
Israel J. Math. \textbf{135} (2003), 83--91. 

\bibitem{Ber}
E. Berger, 
\textit{Hurwitz equivalence in dihedral groups},
Electron. J. Combin. \textbf{18} (2011), no. 1, Paper 45, 16 pp. 

\bibitem{CKS}
J. S. Carter, S. Kamada, and M. Saito, 
\textit{Surfaces in $4$-space}, 
Encyclopaedia of Mathematical Sciences \textbf{142}, 
Low-Dimensional Topology III, Springer, Berlin, 2004. 

\bibitem{Cla}
F. J.-B. J. Clauwens, 
\textit{The adjoint group of an Alexander quandle}, 
arXiv:1011.1587 [math.GR]. 

\bibitem{FRR}
R. Fenn, R. Rim\'{a}nyi, and C. Rourke, 
\textit{The braid-permutation group}, 
Topology \textbf{36} (1997), no.~1, 123--135. 

\bibitem{FR}
R. Fenn and C. Rourke, 
\textit{Racks and links in codimension two}, 
J. Knot Theory Ramifications~\textbf{1} (1992), no. 4, 343--406.

\bibitem{Joy} 
D. Joyce, 
\textit{A classifying invariant of knots, the knot quandle}, 
J. Pure Appl. Algebra \textbf{23} (1982), no. 1, 37--65.

\bibitem{Kam96}
S. Kamada, 
\textit{On braid monodromies of non-simple braided surfaces}, 
Math. Proc. Cambridge Philos. Soc. \textbf{120} (1996), no. 2, 237--245.

\bibitem{Kam02}
S. Kamada, 
\textit{Braid and knot theory in dimension four}, 
Math. Surveys Monogr., \textbf{95}, 
American Mathematical Society, Providence, RI, 2002.

\bibitem{Kam07}
S. Kamada, 
\textit{Braid presentation of virtual knots and welded knots},
Osaka J. Math. \textbf{44} (2007), no. 2, 441--458. 

\bibitem{Kam17}
S. Kamada, 
\textit{Surface-knots in $4$-space}, 
Springer Monographs in Mathematics, Springer, Singapore, 2017. 

\bibitem{KL}
L. H. Kauffman and S. Lambropoulou, 
\textit{Virtual braids}, 
Fund. Math. \textbf{184} (2004), 159--186. 

\bibitem{KT}
V.~S. Kulikov and M. Teicher, 
\textit{Braid monodromy factorizations and diffeomorphism types}, 
Izv. Math. \textbf{64} (2000), no.~2, 311--341; translated from Izv. Ross. Akad. Nauk Ser. Mat. \textbf{64} (2000), no.~2, 89--120. 

\bibitem{Mat} 
S. V. Matveev, 
\textit{Distributive groupoids in knot theory} (Russian), 
Mat. Sb. (N.S.) \textbf{119}(\textbf{161}) (1982), no. 1, 78--88, 160. 

\bibitem{NNSW}
T. Nakamura, Y. Nakanishi, S. Satoh, and K. Wada, 
\textit{Hurwitz equivalence in the universal dihedral quandle}, 
arXiv:2410.11599 [math.GT]. 

\bibitem{Sia} 
C. Sia, 
\textit{Hurwitz equivalence in tuples of dihedral groups, 
dicyclic groups, and semidihedral groups}, 
Electron. J. Combin. \textbf{16} (2009), no. 1, Research Paper 95, 17 pp.
\end{thebibliography}
\end{document}